\documentclass{gtart}

\def\ifplaintex{\expandafter\ifx\csname documentclass\endcsname\relax}

\def\gtp{{\mathsurround=0pt\it $\cal G\mskip-2mu$eometry \&\ 
$\cal T\!\!$opology $\cal P\!$ublications}}  

\def\recd{{\small Received:\qua\receiveddate\ifx\reviseddate\relax
\else\qquad Revised:\qua\reviseddate\fi\par}} 

\def\lognumber#1{\def\thelognumber{#1}}
\def\volumenumber#1{\def\thevolumenumber{#1}}
\def\volumeyear#1{\def\thevolumeyear{#1}}
\def\papernumber#1{\def\thepapernumber{#1}}
\def\pagenumbers#1#2{\def\startpage{#1}\def\finishpage{#2}}
\def\published#1{\def\publishdate{#1}}

\def\received#1{\def\receiveddate{#1}}

\def\accepted#1{\def\accepteddate{#1}}

\def\asciiaddress#1{\def\theasciiaddress{#1}}

\long\def\asciiabstract#1{\long\def\theasciiabstract{#1}}

\let\\\par\let\thelognumber\relax\let\thevolumenumber\relax
\let\thepapernumber\relax\let\thevolumeyear\relax\let\startpage\relax
\let\finishpage\relax\let\publishdate\relax\let\receiveddate\relax
\let\reviseddate\relax\let\accepteddate\relax\let\theasciititle\relax
\let\theasciiauthors\relax\let\theasciiaddress\relax
\let\theasciiabstract\relax

\let\theasciiemail\relax

\ifplaintex
\font\logobig=cmssbx10 scaled 3836
\font\logomed=cmssbx10 scaled 2557
\else
\font\logobig=cmssbx10 scaled 4200
\font\logomed=cmssbx10 scaled 2800
\fi

\long\def\makeagttitle{   
\count0=\startpage
\agt\hfill      
\hbox to 45truept{\vbox to 0pt{\vglue -13truept{\logomed A\kern -.37em{\logobig 
T}\kern -.38em G}\vss}\hss}
\break
{\small Volume \thevolumenumber\ (\thevolumeyear)
\startpage--\finishpage\nl
Published: \publishdate}

\vglue .25truein

{\parskip=0pt\leftskip 0pt plus
1fil\def\\{\par\smallskip}{\Large\bf\thetitle}\par\medskip} \vglue
0.05truein

{\parskip=0pt\leftskip 0pt plus 1fil\def\\{\par}{\sc\theauthors}
\par\medskip}%
 
\vglue 0.03truein 

{\small\leftskip 25truept\rightskip 25truept{\bf Abstract}\stdspace\theabstract

{\bf AMS Classification}\stdspace\theprimaryclass
\ifx\thesecondaryclass\relax\else; \thesecondaryclass\fi\par
{\bf Keywords}\stdspace \thekeywords\par}\vglue 7truept

}   

\ifplaintex
\font\phead=cmsl9 scaled 950
\font\pnum=cmbx10 scaled 913
\font\pfoot=cmsl9 scaled 950
\headline{\vbox to 0pt{\vskip -4.5mm\line{\small\phead\ifnum
\count0=\startpage ISSN 1472-2739 (on-line) 1472-2747 (printed)
\hfill {\pnum\folio}\else\ifodd\count0\def\\{ }%
\ifx\theshorttitle\relax\thetitle\else\theshorttitle\fi\hfill{\pnum\folio}
\else\def\\{ and }{\pnum\folio}\hfill\ifx\theshortauthors\relax\theauthors
\else\theshortauthors\fi\fi\fi}\vss}}
\footline{\vbox to 0pt{\vglue 0mm\line{\small\pfoot\ifnum\count0=\startpage
\copyright\ \gtp\hfill\else
\agt, Volume \thevolumenumber\ (\thevolumeyear)\hfill\fi}\vss}}
\else
\font=cmsl9 scaled 1050
\font\lnum=cmbx10 
\font=cmsl9 scaled 1050
\makeatletter
\def\@oddhead{{\small\lhead\ifnum\count0=\startpage ISSN 1472-2739 
(on-line) 1472-2747 (printed)\hfill {\lnum\number\count0}\else\ifodd\count0
\def\\{ }\ifx\theshorttitle\relax \thetitle \else\theshorttitle\fi\hfill
{\lnum\number\count0}\else\def\\{ and }{\lnum\number\count0}
\hfill\ifx\theshortauthors\relax 
\theauthors\else\theshortauthors\fi\fi\fi}}\def\@evenhead{\@oddhead}
\def\@oddfoot{\small\lfoot\ifnum\count0=\startpage\copyright\ \gtp\hfill\else
\agt, Volume \thevolumenumber\ (\thevolumeyear)\hfill\fi}
\def\@evenfoot{\@oddfoot}
\makeatother
\fi
\let\maketitlepage\makeagttitle
\let\makeshorttitle\maketitlepage
\let\maketitle\maketitlepage

\newwrite\gtoutfile
\long\gdef\makeheadfile{  
{\def\\{, }\def\s{ }
\immediate\openout\gtoutfile head.xxx
\immediate\write\gtoutfile{To: math@arxiv.org}
\immediate\write\gtoutfile{Subject: put OR rep NNNNN:ppppp}
\immediate\write\gtoutfile{--text follows this line--}
\immediate\write\gtoutfile{Proxy-for: \ifx\theasciiauthors\relax
\theauthors\else\theasciiauthors\fi\s<\ifx\theasciiemail\relax\theemail\else\theasciiemail\fi>}
\immediate\write\gtoutfile{\noexpand\\}
\immediate\write\gtoutfile{Authors: \ifx\theasciiauthors\relax
\theauthors\else\theasciiauthors\fi}
{\def\\{ }\immediate\write\gtoutfile{Title: \ifx\theasciititle\relax
\thetitle\else\theasciititle\fi}}
\immediate\write\gtoutfile{Subj-class: GT or SG, GR etc}
\immediate\write\gtoutfile{MSC-class: \theprimaryclass\ifx\thesecondaryclass\relax\else, \thesecondaryclass\fi}
\immediate\write\gtoutfile{Journal-ref: Algebr. Geom. Topol. \thevolumenumber\s
(\thevolumeyear) \startpage-\finishpage}
\immediate\write\gtoutfile{Comments: Published by Algebraic and
Geometric Topology at}
\immediate\write\gtoutfile{\s\s\s  http://www.maths.warwick.ac.uk/agt/AGTVol\thevolumenumber/agt-\thevolumenumber-\thepapernumber.abs.html}
\immediate\write\gtoutfile{\noexpand\\}
\immediate\write\gtoutfile{}
\ifx\theasciiabstract\relax
\immediate\write\gtoutfile{\theabstract}\else
\immediate\write\gtoutfile{\theasciiabstract}\fi
\immediate\write\gtoutfile{}
\immediate\write\gtoutfile{\noexpand\\}
\immediate\write\gtoutfile{}
\immediate\closeout\gtoutfile}}  

\def\maketitlepage{\makeagttitle\makeheadfile}
\let\makeshorttitle\maketitlepage
\let\maketitle\maketitlepage

\def\ifplaintex{\expandafter\ifx\csname documentclass\endcsname\relax}

\def\gtp{{\mathsurround=0pt\it $\cal G\mskip-2mu$eometry \&\ 
$\cal T\!\!$opology $\cal P\!$ublications}}  

\def\recd{{\small Received:\qua\receiveddate\ifx\reviseddate\relax
\else\qquad Revised:\qua\reviseddate\fi\par}} 

\def\lognumber#1{\def\thelognumber{#1}}
\def\volumenumber#1{\def\thevolumenumber{#1}}
\def\volumeyear#1{\def\thevolumeyear{#1}}
\def\papernumber#1{\def\thepapernumber{#1}}
\def\pagenumbers#1#2{\def\startpage{#1}\def\finishpage{#2}}
\def\published#1{\def\publishdate{#1}}

\def\received#1{\def\receiveddate{#1}}

\def\accepted#1{\def\accepteddate{#1}}

\def\asciiaddress#1{\def\theasciiaddress{#1}}

\long\def\asciiabstract#1{\long\def\theasciiabstract{#1}}

\let\\\par\let\thelognumber\relax\let\thevolumenumber\relax
\let\thepapernumber\relax\let\thevolumeyear\relax\let\startpage\relax
\let\finishpage\relax\let\publishdate\relax\let\receiveddate\relax
\let\reviseddate\relax\let\accepteddate\relax\let\theasciititle\relax
\let\theasciiauthors\relax\let\theasciiaddress\relax
\let\theasciiabstract\relax

\let\theasciiemail\relax

\ifplaintex
\font\logobig=cmssbx10 scaled 3836
\font\logomed=cmssbx10 scaled 2557
\else
\font\logobig=cmssbx10 scaled 4200
\font\logomed=cmssbx10 scaled 2800
\fi

\long\def\makeagttitle{   
\count0=\startpage
\agt\hfill      
\hbox to 45truept{\vbox to 0pt{\vglue -13truept{\logomed A\kern -.37em{\logobig 
T}\kern -.38em G}\vss}\hss}
\break
{\small Volume \thevolumenumber\ (\thevolumeyear)
\startpage--\finishpage\nl
Published: \publishdate}

\vglue .25truein

{\parskip=0pt\leftskip 0pt plus
1fil\def\\{\par\smallskip}{\Large\bf\thetitle}\par\medskip} \vglue
0.05truein

{\parskip=0pt\leftskip 0pt plus 1fil\def\\{\par}{\sc\theauthors}
\par\medskip}%
 
\vglue 0.03truein 

{\small\leftskip 25truept\rightskip 25truept{\bf Abstract}\stdspace\theabstract

{\bf AMS Classification}\stdspace\theprimaryclass
\ifx\thesecondaryclass\relax\else; \thesecondaryclass\fi\par
{\bf Keywords}\stdspace \thekeywords\par}\vglue 7truept

}   

\ifplaintex
\font\phead=cmsl9 scaled 950
\font\pnum=cmbx10 scaled 913
\font\pfoot=cmsl9 scaled 950
\headline{\vbox to 0pt{\vskip -4.5mm\line{\small\phead\ifnum
\count0=\startpage ISSN 1472-2739 (on-line) 1472-2747 (printed)
\hfill {\pnum\folio}\else\ifodd\count0\def\\{ }%
\ifx\theshorttitle\relax\thetitle\else\theshorttitle\fi\hfill{\pnum\folio}
\else\def\\{ and }{\pnum\folio}\hfill\ifx\theshortauthors\relax\theauthors
\else\theshortauthors\fi\fi\fi}\vss}}
\footline{\vbox to 0pt{\vglue 0mm\line{\small\pfoot\ifnum\count0=\startpage
\copyright\ \gtp\hfill\else
\agt, Volume \thevolumenumber\ (\thevolumeyear)\hfill\fi}\vss}}
\else
\font=cmsl9 scaled 1050
\font\lnum=cmbx10 
\font=cmsl9 scaled 1050
\makeatletter
\def\@oddhead{{\small\lhead\ifnum\count0=\startpage ISSN 1472-2739 
(on-line) 1472-2747 (printed)\hfill {\lnum\number\count0}\else\ifodd\count0
\def\\{ }\ifx\theshorttitle\relax \thetitle \else\theshorttitle\fi\hfill
{\lnum\number\count0}\else\def\\{ and }{\lnum\number\count0}
\hfill\ifx\theshortauthors\relax 
\theauthors\else\theshortauthors\fi\fi\fi}}\def\@evenhead{\@oddhead}
\def\@oddfoot{\small\lfoot\ifnum\count0=\startpage\copyright\ \gtp\hfill\else
\agt, Volume \thevolumenumber\ (\thevolumeyear)\hfill\fi}
\def\@evenfoot{\@oddfoot}
\makeatother
\fi
\let\maketitlepage\makeagttitle
\let\makeshorttitle\maketitlepage
\let\maketitle\maketitlepage

\newwrite\gtoutfile
\long\gdef\makeheadfile{  
{\def\\{, }\def\s{ }
\immediate\openout\gtoutfile head.xxx
\immediate\write\gtoutfile{To: math@arxiv.org}
\immediate\write\gtoutfile{Subject: put OR rep NNNNN:ppppp}
\immediate\write\gtoutfile{--text follows this line--}
\immediate\write\gtoutfile{Proxy-for: \ifx\theasciiauthors\relax
\theauthors\else\theasciiauthors\fi\s<\ifx\theasciiemail\relax\theemail\else\theasciiemail\fi>}
\immediate\write\gtoutfile{\noexpand\\}
\immediate\write\gtoutfile{Authors: \ifx\theasciiauthors\relax
\theauthors\else\theasciiauthors\fi}
{\def\\{ }\immediate\write\gtoutfile{Title: \ifx\theasciititle\relax
\thetitle\else\theasciititle\fi}}
\immediate\write\gtoutfile{Subj-class: GT or SG, GR etc}
\immediate\write\gtoutfile{MSC-class: \theprimaryclass\ifx\thesecondaryclass\relax\else, \thesecondaryclass\fi}
\immediate\write\gtoutfile{Journal-ref: Algebr. Geom. Topol. \thevolumenumber\s
(\thevolumeyear) \startpage-\finishpage}
\immediate\write\gtoutfile{Comments: Published by Algebraic and
Geometric Topology at}
\immediate\write\gtoutfile{\s\s\s  http://www.maths.warwick.ac.uk/agt/AGTVol\thevolumenumber/agt-\thevolumenumber-\thepapernumber.abs.html}
\immediate\write\gtoutfile{\noexpand\\}
\immediate\write\gtoutfile{}
\ifx\theasciiabstract\relax
\immediate\write\gtoutfile{\theabstract}\else
\immediate\write\gtoutfile{\theasciiabstract}\fi
\immediate\write\gtoutfile{}
\immediate\write\gtoutfile{\noexpand\\}
\immediate\write\gtoutfile{}
\immediate\closeout\gtoutfile}}  

\def\maketitlepage{\makeagttitle\makeheadfile}
\let\makeshorttitle\maketitlepage
\let\maketitle\maketitlepage

\def\ifplaintex{\expandafter\ifx\csname documentclass\endcsname\relax}

\def\gtp{{\mathsurround=0pt\it $\cal G\mskip-2mu$eometry \&\ 
$\cal T\!\!$opology $\cal P\!$ublications}}  

\def\recd{{\small Received:\qua\receiveddate\ifx\reviseddate\relax
\else\qquad Revised:\qua\reviseddate\fi\par}} 

\def\lognumber#1{\def\thelognumber{#1}}
\def\volumenumber#1{\def\thevolumenumber{#1}}
\def\volumeyear#1{\def\thevolumeyear{#1}}
\def\papernumber#1{\def\thepapernumber{#1}}
\def\pagenumbers#1#2{\def\startpage{#1}\def\finishpage{#2}}
\def\published#1{\def\publishdate{#1}}

\def\received#1{\def\receiveddate{#1}}

\def\accepted#1{\def\accepteddate{#1}}

\def\asciiaddress#1{\def\theasciiaddress{#1}}

\long\def\asciiabstract#1{\long\def\theasciiabstract{#1}}

\let\\\par\let\thelognumber\relax\let\thevolumenumber\relax
\let\thepapernumber\relax\let\thevolumeyear\relax\let\startpage\relax
\let\finishpage\relax\let\publishdate\relax\let\receiveddate\relax
\let\reviseddate\relax\let\accepteddate\relax\let\theasciititle\relax
\let\theasciiauthors\relax\let\theasciiaddress\relax
\let\theasciiabstract\relax

\let\theasciiemail\relax

\ifplaintex
\font\logobig=cmssbx10 scaled 3836
\font\logomed=cmssbx10 scaled 2557
\else
\font\logobig=cmssbx10 scaled 4200
\font\logomed=cmssbx10 scaled 2800
\fi

\long\def\makeagttitle{   
\count0=\startpage
\agt\hfill      
\hbox to 45truept{\vbox to 0pt{\vglue -13truept{\logomed A\kern -.37em{\logobig 
T}\kern -.38em G}\vss}\hss}
\break
{\small Volume \thevolumenumber\ (\thevolumeyear)
\startpage--\finishpage\nl
Published: \publishdate}

\vglue .25truein

{\parskip=0pt\leftskip 0pt plus
1fil\def\\{\par\smallskip}{\Large\bf\thetitle}\par\medskip} \vglue
0.05truein

{\parskip=0pt\leftskip 0pt plus 1fil\def\\{\par}{\sc\theauthors}
\par\medskip}%
 
\vglue 0.03truein 

{\small\leftskip 25truept\rightskip 25truept{\bf Abstract}\stdspace\theabstract

{\bf AMS Classification}\stdspace\theprimaryclass
\ifx\thesecondaryclass\relax\else; \thesecondaryclass\fi\par
{\bf Keywords}\stdspace \thekeywords\par}\vglue 7truept

}   

\ifplaintex
\font\phead=cmsl9 scaled 950
\font\pnum=cmbx10 scaled 913
\font\pfoot=cmsl9 scaled 950
\headline{\vbox to 0pt{\vskip -4.5mm\line{\small\phead\ifnum
\count0=\startpage ISSN 1472-2739 (on-line) 1472-2747 (printed)
\hfill {\pnum\folio}\else\ifodd\count0\def\\{ }%
\ifx\theshorttitle\relax\thetitle\else\theshorttitle\fi\hfill{\pnum\folio}
\else\def\\{ and }{\pnum\folio}\hfill\ifx\theshortauthors\relax\theauthors
\else\theshortauthors\fi\fi\fi}\vss}}
\footline{\vbox to 0pt{\vglue 0mm\line{\small\pfoot\ifnum\count0=\startpage
\copyright\ \gtp\hfill\else
\agt, Volume \thevolumenumber\ (\thevolumeyear)\hfill\fi}\vss}}
\else
\font=cmsl9 scaled 1050
\font\lnum=cmbx10 
\font=cmsl9 scaled 1050
\makeatletter
\def\@oddhead{{\small\lhead\ifnum\count0=\startpage ISSN 1472-2739 
(on-line) 1472-2747 (printed)\hfill {\lnum\number\count0}\else\ifodd\count0
\def\\{ }\ifx\theshorttitle\relax \thetitle \else\theshorttitle\fi\hfill
{\lnum\number\count0}\else\def\\{ and }{\lnum\number\count0}
\hfill\ifx\theshortauthors\relax 
\theauthors\else\theshortauthors\fi\fi\fi}}\def\@evenhead{\@oddhead}
\def\@oddfoot{\small\lfoot\ifnum\count0=\startpage\copyright\ \gtp\hfill\else
\agt, Volume \thevolumenumber\ (\thevolumeyear)\hfill\fi}
\def\@evenfoot{\@oddfoot}
\makeatother
\fi
\let\maketitlepage\makeagttitle
\let\makeshorttitle\maketitlepage
\let\maketitle\maketitlepage

\newwrite\gtoutfile
\long\gdef\makeheadfile{  
{\def\\{, }\def\s{ }
\immediate\openout\gtoutfile head.xxx
\immediate\write\gtoutfile{To: math@arxiv.org}
\immediate\write\gtoutfile{Subject: put OR rep NNNNN:ppppp}
\immediate\write\gtoutfile{--text follows this line--}
\immediate\write\gtoutfile{Proxy-for: \ifx\theasciiauthors\relax
\theauthors\else\theasciiauthors\fi\s<\ifx\theasciiemail\relax\theemail\else\theasciiemail\fi>}
\immediate\write\gtoutfile{\noexpand\\}
\immediate\write\gtoutfile{Authors: \ifx\theasciiauthors\relax
\theauthors\else\theasciiauthors\fi}
{\def\\{ }\immediate\write\gtoutfile{Title: \ifx\theasciititle\relax
\thetitle\else\theasciititle\fi}}
\immediate\write\gtoutfile{Subj-class: GT or SG, GR etc}
\immediate\write\gtoutfile{MSC-class: \theprimaryclass\ifx\thesecondaryclass\relax\else, \thesecondaryclass\fi}
\immediate\write\gtoutfile{Journal-ref: Algebr. Geom. Topol. \thevolumenumber\s
(\thevolumeyear) \startpage-\finishpage}
\immediate\write\gtoutfile{Comments: Published by Algebraic and
Geometric Topology at}
\immediate\write\gtoutfile{\s\s\s  http://www.maths.warwick.ac.uk/agt/AGTVol\thevolumenumber/agt-\thevolumenumber-\thepapernumber.abs.html}
\immediate\write\gtoutfile{\noexpand\\}
\immediate\write\gtoutfile{}
\ifx\theasciiabstract\relax
\immediate\write\gtoutfile{\theabstract}\else
\immediate\write\gtoutfile{\theasciiabstract}\fi
\immediate\write\gtoutfile{}
\immediate\write\gtoutfile{\noexpand\\}
\immediate\write\gtoutfile{}
\immediate\closeout\gtoutfile}}  

\def\maketitlepage{\makeagttitle\makeheadfile}
\let\makeshorttitle\maketitlepage
\let\maketitle\maketitlepage

\def\ifplaintex{\expandafter\ifx\csname documentclass\endcsname\relax}

\def\gtp{{\mathsurround=0pt\it $\cal G\mskip-2mu$eometry \&\ 
$\cal T\!\!$opology $\cal P\!$ublications}}  

\def\recd{{\small Received:\qua\receiveddate\ifx\reviseddate\relax
\else\qquad Revised:\qua\reviseddate\fi\par}} 

\def\lognumber#1{\def\thelognumber{#1}}
\def\volumenumber#1{\def\thevolumenumber{#1}}
\def\volumeyear#1{\def\thevolumeyear{#1}}
\def\papernumber#1{\def\thepapernumber{#1}}
\def\pagenumbers#1#2{\def\startpage{#1}\def\finishpage{#2}}
\def\published#1{\def\publishdate{#1}}

\def\received#1{\def\receiveddate{#1}}

\def\accepted#1{\def\accepteddate{#1}}

\def\asciiaddress#1{\def\theasciiaddress{#1}}

\long\def\asciiabstract#1{\long\def\theasciiabstract{#1}}

\let\\\par\let\thelognumber\relax\let\thevolumenumber\relax
\let\thepapernumber\relax\let\thevolumeyear\relax\let\startpage\relax
\let\finishpage\relax\let\publishdate\relax\let\receiveddate\relax
\let\reviseddate\relax\let\accepteddate\relax\let\theasciititle\relax
\let\theasciiauthors\relax\let\theasciiaddress\relax
\let\theasciiabstract\relax

\let\theasciiemail\relax

\ifplaintex
\font\logobig=cmssbx10 scaled 3836
\font\logomed=cmssbx10 scaled 2557
\else
\font\logobig=cmssbx10 scaled 4200
\font\logomed=cmssbx10 scaled 2800
\fi

\long\def\makeagttitle{   
\count0=\startpage
\agt\hfill      
\hbox to 45truept{\vbox to 0pt{\vglue -13truept{\logomed A\kern -.37em{\logobig 
T}\kern -.38em G}\vss}\hss}
\break
{\small Volume \thevolumenumber\ (\thevolumeyear)
\startpage--\finishpage\nl
Published: \publishdate}

\vglue .25truein

{\parskip=0pt\leftskip 0pt plus
1fil\def\\{\par\smallskip}{\Large\bf\thetitle}\par\medskip} \vglue
0.05truein

{\parskip=0pt\leftskip 0pt plus 1fil\def\\{\par}{\sc\theauthors}
\par\medskip}%
 
\vglue 0.03truein 

{\small\leftskip 25truept\rightskip 25truept{\bf Abstract}\stdspace\theabstract

{\bf AMS Classification}\stdspace\theprimaryclass
\ifx\thesecondaryclass\relax\else; \thesecondaryclass\fi\par
{\bf Keywords}\stdspace \thekeywords\par}\vglue 7truept

}   

\ifplaintex
\font\phead=cmsl9 scaled 950
\font\pnum=cmbx10 scaled 913
\font\pfoot=cmsl9 scaled 950
\headline{\vbox to 0pt{\vskip -4.5mm\line{\small\phead\ifnum
\count0=\startpage ISSN 1472-2739 (on-line) 1472-2747 (printed)
\hfill {\pnum\folio}\else\ifodd\count0\def\\{ }%
\ifx\theshorttitle\relax\thetitle\else\theshorttitle\fi\hfill{\pnum\folio}
\else\def\\{ and }{\pnum\folio}\hfill\ifx\theshortauthors\relax\theauthors
\else\theshortauthors\fi\fi\fi}\vss}}
\footline{\vbox to 0pt{\vglue 0mm\line{\small\pfoot\ifnum\count0=\startpage
\copyright\ \gtp\hfill\else
\agt, Volume \thevolumenumber\ (\thevolumeyear)\hfill\fi}\vss}}
\else
\font=cmsl9 scaled 1050
\font\lnum=cmbx10 
\font=cmsl9 scaled 1050
\makeatletter
\def\@oddhead{{\small\lhead\ifnum\count0=\startpage ISSN 1472-2739 
(on-line) 1472-2747 (printed)\hfill {\lnum\number\count0}\else\ifodd\count0
\def\\{ }\ifx\theshorttitle\relax \thetitle \else\theshorttitle\fi\hfill
{\lnum\number\count0}\else\def\\{ and }{\lnum\number\count0}
\hfill\ifx\theshortauthors\relax 
\theauthors\else\theshortauthors\fi\fi\fi}}\def\@evenhead{\@oddhead}
\def\@oddfoot{\small\lfoot\ifnum\count0=\startpage\copyright\ \gtp\hfill\else
\agt, Volume \thevolumenumber\ (\thevolumeyear)\hfill\fi}
\def\@evenfoot{\@oddfoot}
\makeatother
\fi
\let\maketitlepage\makeagttitle
\let\makeshorttitle\maketitlepage
\let\maketitle\maketitlepage

\newwrite\gtoutfile
\long\gdef\makeheadfile{  
{\def\\{, }\def\s{ }
\immediate\openout\gtoutfile head.xxx
\immediate\write\gtoutfile{To: math@arxiv.org}
\immediate\write\gtoutfile{Subject: put OR rep NNNNN:ppppp}
\immediate\write\gtoutfile{--text follows this line--}
\immediate\write\gtoutfile{Proxy-for: \ifx\theasciiauthors\relax
\theauthors\else\theasciiauthors\fi\s<\ifx\theasciiemail\relax\theemail\else\theasciiemail\fi>}
\immediate\write\gtoutfile{\noexpand\\}
\immediate\write\gtoutfile{Authors: \ifx\theasciiauthors\relax
\theauthors\else\theasciiauthors\fi}
{\def\\{ }\immediate\write\gtoutfile{Title: \ifx\theasciititle\relax
\thetitle\else\theasciititle\fi}}
\immediate\write\gtoutfile{Subj-class: GT or SG, GR etc}
\immediate\write\gtoutfile{MSC-class: \theprimaryclass\ifx\thesecondaryclass\relax\else, \thesecondaryclass\fi}
\immediate\write\gtoutfile{Journal-ref: Algebr. Geom. Topol. \thevolumenumber\s
(\thevolumeyear) \startpage-\finishpage}
\immediate\write\gtoutfile{Comments: Published by Algebraic and
Geometric Topology at}
\immediate\write\gtoutfile{\s\s\s  http://www.maths.warwick.ac.uk/agt/AGTVol\thevolumenumber/agt-\thevolumenumber-\thepapernumber.abs.html}
\immediate\write\gtoutfile{\noexpand\\}
\immediate\write\gtoutfile{}
\ifx\theasciiabstract\relax
\immediate\write\gtoutfile{\theabstract}\else
\immediate\write\gtoutfile{\theasciiabstract}\fi
\immediate\write\gtoutfile{}
\immediate\write\gtoutfile{\noexpand\\}
\immediate\write\gtoutfile{}
\immediate\closeout\gtoutfile}}  

\def\maketitlepage{\makeagttitle\makeheadfile}
\let\makeshorttitle\maketitlepage
\let\maketitle\maketitlepage

\lognumber{25}
\volumenumber{2}
\volumeyear{2002}
\papernumber{25}
\published{28 June 2002}
\pagenumbers{519}{536}
\received{11 October 2001}
\accepted{20 June 2002}

\usepackage{amssymb}
\input epsf

\def\ppcm{\mathop{\rm ppcm}}
\def\pgcd{\mathop{\rm pgcd}}
\def\lcm{\mathop{\rm lcm}}

\newtheorem{Lemme}{Lemme}[section]
\newtheorem{The}[Lemme]{Th\'eor\`eme}
\newtheorem{Prop}[Lemme]{Proposition}
\newtheorem{Cor}[Lemme]{Corollaire}

\theoremstyle{definition}
\newtheorem{Def}[Lemme]{D\'efinition}

\newtheorem{Que}{Question}

\begin{document}
\title{Morphismes injectifs entre groupes d'Artin-Tits}
\author{Eddy Godelle}
\address{LAMFA CNRS 2270, Universit\'e de Picardie-Jules
 Verne\\Facult\'e de Math\'ematiques et d'Informatique\\33 rue Saint-Leu,
 80000 Amiens,  France}
\asciiaddress{LAMFA CNRS 2270, Universite de Picardie-Jules
 Verne\\Faculte de Mathematiques et d'Informatique\\33 rue Saint-Leu,
 80000 Amiens,  France}

\email{eddy.godelle@u-picardie.fr}
\url{http://www.mathinfo.u-picardie.fr/godelle}

\begin{abstract}  We construct a family of morphisms between Artin-Tits groups which generalise the ones constructed by J. Crisp in \cite{Cri}. We show that their restrictions to the positive Artin monoids respect normal forms, and that for Artin-Tits groups of type FC, these morphisms are injective. The proof of the second result uses the Deligne Complex, and the normal cube paths constructed in \cite{NiR} and \cite{ChA}.

{\bf R\'esum\'e}\qua On construit une classe de morphismes entre groupes
 d'Artin-Tits qui g\'en\'eralise
celle construite par J. Crisp dans \cite{Cri}. On montre que leurs
 restrictions aux mono\"\i des respectent les formes normales, et que pour les
 groupes d'Artin-Tits de type FC ces morphismes sont injectifs. La
 d\'emonstration du second r\'esultat utilise
 le complexe de Deligne et les chemins cubiques normaux construits dans
 \cite{NiR} et \cite{ChA}. \end{abstract}

\asciiabstract{We construct a family of morphisms between Artin-Tits
groups which generalise the ones constructed by J. Crisp in [Injective
maps between Artin groups, Proceedings of the Special Year in
Geometric Group Theory, Berlin, (1999), 119 -- 138]. We show that
their restrictions to the positive Artin monoids respect normal forms,
and that for Artin-Tits groups of type FC, these morphisms are
injective. The proof of the second result uses the Deligne Complex,
and the normal cube paths constructed in [G. Niblo and L. Reeves, The
geometry of cube complexes and the complexity of their fundamental
groups, Topology 37 (1998) 621-633] and [J.A. Altobelli and
R. Charney, A geometric Rational Form for Artin Groups of FC type,
Geom. Dedicata, 79 (2000) 277-289].

Resume: 

On construit une classe de morphismes entre groupes
d'Artin-Tits qui generalise celle construite par J. Crisp dans
[Injective maps between Artin groups, Proceedings of the Special Year
in Geometric Group Theory, Berlin, (1999), 119 -- 138]. On montre que
leurs restrictions aux monoides respectent les formes normales, et que
pour les groupes d'Artin-Tits de type FC ces morphismes sont
injectifs. La demonstration du second resultat utilise le complexe de
Deligne et les chemins cubiques normaux construits dans [G. Niblo et
L. Reeves, The geometry of cube complexes and the complexity of their
fundamental groups, Topology 37 (1998) 621-633] et [J.A. Altobelli et
R. Charney, A geometric Rational Form for Artin Groups of FC type,
Geom. Dedicata, 79 (2000) 277-289].}

\maketitle

\let\\\par

\section*{Introduction} Soit $S$ un ensemble fini et $M=
(m_{s,t})_{s,t\in S}$ une matrice sym\'etrique \`a coefficients dans
$\mathbb{N}\cup\{\infty\}-\{0\}$ telle que $m_{s,s} = 1$  pour tout $s$
de $S$ et $m_{s,t} \neq 1$ pour $s\not=t$ dans $S$. On note $A_S$ le groupe engendr\'e par $S$ et les relations, dites
``de tresses'', $\underbrace{sts\cdots}_{m_{s,t}\ termes} =
\underbrace{tst\cdots}_{m_{s,t}\ termes}$ pour tout couple $(s,t)$
d'\'el\'ements distincts de $S$ tels que $m_{s,t}\not=\infty$: $$A_S =
\langle S|\underbrace{sts\cdots}_{m_{s,t}\ termes} =
\underbrace{tst\cdots}_{m_{s,t}\ termes}\ ;\ \forall s,t\in S, s\not= t\
et\ m_{s,t}\not=\infty \rangle.$$ La paire $(A_S,S)$ s'appelle un syst\`eme
d'Artin-Tits et $A_S$ un groupe d'Artin-Tits (relativement \`a $S$). On note $A^+_S$ le sous-mono\"\i de de $A_S$
engendr\'e par $S$. Ce mono\"\i de $A_S^+$ poss\`ede la m\^eme
pr\'esentation que $A_S$ mais en tant que mono\"\i de
(\cite{Par4}). Puisque les relations de tresses sont homog\`enes, $A_S^+$
est naturellement muni d'une fonction longueur $\ell$ compatible au
produit. On appelle graphe de $S$ et $M$, not\'e $\Gamma_S$, le graphe
\'etiquet\'e dont l'ensemble des sommets est $S$ et dont les ar\^etes
sont les paires $\{s,t\}$ d'\'el\'ements distincts de $S$ telles que
$m_{s,t}\not= 2$, que l'on \'etiquette par $m_{s,t}$.\\ On appelle
sous-groupe parabolique standard tout sous-groupe de $A_S$
engendr\'e par une partie $T$ de $S$; on note $A_T$ un tel
sous-groupe.  Van der Lek a prouv\'e dans \cite{VdL} que pour toute
partie $T$ de $S$, la paire
$(A_T,T)$ est un syst\`eme d'Artin-Tits pour la matrice $(m_{s,t})_{s,t\in T}$.
Un sous-groupe parabolique de $A_S$ est un
sous-groupe de $A_S$ conjugu\'e \`a un sous-groupe parabolique
standard de $A_S$.\\ Lorsque ce graphe
est connexe, on dit que $S$ est ind\'ecomposable. Lorsque l'on ajoute
\`a la pr\'esentation de $A_S$ les relations $s^2 = 1$ pour $s\in S$, on
obtient un groupe de Coxeter $W_S$. On dit que $A_S$ (ou $S$ par abus) est
de type sph\'erique lorsque son groupe de Coxeter associ\'e est fini.\\
\noindent Si $u$ et $v$ sont dans $A_S^+$, on notera $u\prec v$ (resp.\
$v\succ u$) pour dire que $u$ divise $v$ \`a gauche (resp.\ \`a
droite); on notera $u\land_\prec v$ (resp.\
$v\land_\succ u$) leur pgcd  relativement \`a
$\prec $ (resp.\
$\succ$) et $u\lor_\prec v$ (resp.\
$v\lor_\succ u$) leur ppcm relativement \`a
$\prec $ (resp.\
$\succ$) s'il existe; s'il n'existe pas, on pose $u\lor_\prec v =
\infty$ (resp.\
$v\lor_\succ u = \infty$). Pour $m\in \mathbb{N}$, on d\'esignera par $[u,v\rangle^m$ le produit $\underbrace{uvu\cdots}_{m\ termes}$. Ainsi les relations de tresses s'\'ecrivent $[s,t\rangle^{m_{s,t}} = [t,s\rangle^{m_{s,t}}$.\\ Rappelons le r\'esultat classique suivant d\'emontr\'e par Brieskorn et Saito dans \cite{BrS}~: un groupe
 d'Artin-Tits $A_S$ est de type sph\'erique si et seulement si le
 $\ppcm$ de $S$ \`a gauche existe; dans ce cas, le ppcm \`a droite existe
 aussi et est \'egal  au $\ppcm$ \`a gauche. On note $\Delta_S$ ou
 simplement $\Delta$ cet \'el\'ement.   
\begin{Def}\label{defLcm}
Soit $A_S$ et $A_{S'}$ deux groupes d'Artin-Tits et soit $p$ une
 application de $S$ dans $ \mathcal{P}(S')-\{\emptyset\}$, les parties
 non vides de $S'$, telle
 que \\ {\bf(L0)}\qua si $s\not=t\in S$ alors $p(s)$ et $p(t)$ sont
 disjointes;\\ {\bf(L1)}\qua pour $s\in S$, $p(s)$ est de type sph\'erique;\\
{\bf (L2)}\qua  si $s\not=t\in S$ avec $m_{s,t}\not= \infty$, on a\\
 \begin{center} $[\Delta_{p(s)},\Delta_{p(t)}\rangle^{m_{s,t}} = [\Delta_{p(t)},\Delta_{p(s)}\rangle^{m_{s,t}} = \Delta_{p(s)}\lor_\prec \Delta_{p(t)}$ dans $A_{S'}^+$\end{center} {\bf(L3)}\qua si $s\not=t\in S$ avec $m_{s,t}= \infty$, alors \begin{center}$\left\{\begin {array}{l} \forall u\in p(s),\ \{u\}\cup p(t)\textrm{ n'est pas de type sph\'erique},\\ \forall u\in p(t),\ \{u\}\cup p(s)\textrm{ n'est pas de type sph\'erique.}\end{array}\right.$ \end{center} On peut alors d\'efinir homomorphisme $\varphi_p$ de $A_S$ dans $A_{S'}$ par $\varphi_p(s)= \Delta_{p(s)}$ pour $s\in S$. Un morphisme provenant d'une telle construction sera appel\'e un LCM-homomorphisme.      
\end{Def}
Par rapport \`a la d\'efinition 2.1 de \cite{Cri}, on a
 ajout\'e la condition {\bf(L3)} qui autorise les liaisons
 infinies. Notons aussi que la d\'efinition 2.1 de \cite{Cri} est donn\'ee dans le cadre des mono\"\i des (voir la definition \ref{deflcmcri} ci-dessous) ce qui est \'equivalent. La proposition suivante g\'en\'eralise la proposition 2.3 de
 \cite{Cri}; voir \'egalement le th\'eor\`eme 14 de \cite{Cri2} pour le
 cas ``sym\'etrique''.
\begin{Prop}\label{propdeb1}
Soit $\varphi_p: A_S\to A_{S'}$ un LCM-homomorphisme. Alors $\varphi_p$ induit un
 homomorphisme injectif $\varphi_{p,W}: W_S\to W_{S'}$.
\end{Prop}
Ce r\'esultat nous permet, comme dans \cite{Cri}, de voir les LCM-homomorphismes
comme des applications entre groupes fondamentaux induites par des applications
simpliciales injectives  (\textit{cf.}
proposition \ref{prorep}).\\
\begin{Def}\label{deffc}\cite{ChD2}\qua On dit qu'un groupe d'Artin-Tits
 $A_S$ est de type $FC$ si et seulement si la condition ci-dessous est
 v\'erifi\'ee pour toute partie $T$ de $S$:\begin{center}
               $\forall s,t\in T$, $m_{s,t}\not=\infty \Rightarrow A_T$ est de type  sph\'erique.\end{center}
\end{Def}
\noindent L'un des deux r\'esultats principaux de cette note est le suivant:   
\begin{The}\label{Thdeb2} Soit $(A_S,S)$ et $(A_{S'},S')$ deux syst\`emes d'Artin-Tits de type FC
et $\varphi_p: A_S\to A_{S'}$ un LCM-homomorphisme. Alors $\varphi_p$ est injective.
\end{The}
Dans \cite{Cri2}, Crisp prouve que le sous-groupe des points fixes d'un groupe d'Artin-Tits de type FC sous l'action d'un groupe de sym\'etries de son graphe est  aussi un groupe d'Artin-Tits de type FC. Il est assez facile de voir que le morphisme construit par Crisp pour la preuve de son r\'esultat v\'erifie les axiomes de la d\'efinition \ref{defLcm} et est donc un LCM-homomorphisme. Le th\'eor\`eme \ref{Thdeb2} peut donc \^etre vu comme une g\'en\'eralisation du r\'esultat de Crisp relatif \`a l'injectivit\'e. \\ 
\begin{Lemme}\label{lemsecpi} Soit $(A_S,S)$ un syst\`eme d'Artin-Tits, soit $W_S$ son groupe de
 Coxeter associ\'e et $i:A_S\to W_S $ la surjection canonique. Si $w$
 est dans $W_S$, alors $i^{-1}(w)\cap A_S^+$ poss\`ede un unique repr\'esentant $Sec(w)$ de
 longueur minimale. On peut ainsi d\'efinir une  section
ensembliste de $i$, not\'ee $Sec$ dont l'image not\'ee $A_{S,red}$ est
 dans $A_S^+$. Ses \'el\'ements sont caract\'eris\'es par le fait qu'aucune de leurs
 \'ecritures dans $A_S^+$ ne fait appara\^\i tre de carr\'es d'un
 \'el\'ement de $S$. 
\end{Lemme}
Cette section ensembliste est construite gr\^ace au lemme d'\'echange (voir \cite{Bou} chapitre 4). Lorsque $A_S$ est de type sph\'erique, alors $\Delta_S$ est l'image par cette section de l'\'el\'ement de plus grande longueur de $W_S$. Les \'el\'ements de $A_{S,red}$ seront dit r\'eduits ou encore minimaux. La derni\`ere
assertion du lemme implique en particulier que $A_{S,red}$ est stable par
division \`a gauche et \`a droite.
\begin{Lemme}\label{lemalpha}{\rm\cite{BrS,Del,Mic}}\qua Soit $(A_S,S)$ un syst\`eme d'Artin-Tits. Pour tout
 \'el\'ement $g$ de $A_S^+$, l'ensemble $\{h\in A_{S,red}| h\prec g\}$
 poss\`ede un plus grand \'el\'ement $\alpha(g)$ pour la division \`a
 gauche. De plus pour $g_1,g_2$ dans $A_S^+$, on a $\alpha(g_1g_2) = \alpha(g_1\alpha(g_2))$.
\end{Lemme}
Cette fonction $\alpha$ permet de d\'efinir une forme normale sur $A_S^+$:
on dira que la suite $(g_1,\cdots, g_n)$ est la forme normale de $g\in
A_S^+$ si $g = g_1\cdots g_n$ o\'u aucun $g_i$ ne vaut $1$ et $g_i = \alpha(g_i\cdots g_n)$ pour
tout $i$; cette
d\'ecomposition est unique d'apr\`es le lemme ci-dessus.\\

\noindent Si $\varphi_p: A_S\to A_{S'}$ est un LCM-homormorphisme alors
l'image par $\varphi$ d'un
\'el\'ement de $A^+_S$ est dans $A^+_{S'}$. Le second r\'esultat que nous allons prouver est le suivant:
\begin{The} \label{cor2deb} Soit $\varphi_p: A_S\to A_{S'}$ un
 LCM-homomorphisme homomorphisme. Alors
 $\varphi_p$ est compatible avec la forme normale: si $(g_1,\cdots,g_n)$ est la forme
 normale de $g\in A_S^+$ alors $(\varphi_p(g_1),\cdots,\varphi_p(g_n))$ est la
 forme normale de $\varphi_p(g)$.
\end{The}
Nous montrerons en fait ce r\'esultat pour une famille un peu plus large de morphismes .\\

Dans la premi\`ere partie nous rappelons les r\'esultats utiles sur les groupes
d'Artin-Tits et nous introduisons la notion de
lcm-homomorphisme; celle-ci g\'en\'eralise celle de LCM-homomorphisme. La seconde partie
est consacr\'ee aux preuves de la proposition \ref{propdeb1} et du th\'eor\`eme
\ref{cor2deb}. Enfin dans la derni\`ere partie, nous introduisons le complexe
de Deligne et nous prouvons le th\'eor\`eme \ref{Thdeb2}.  
\section{G\'en\'eralit\'es}
\subsection{Mono\"\i des d'Artin-Tits} 
\begin{Lemme}{\rm\cite{BrS}}\qua
Soit $(A_S,S)$ un syst\`eme d'Artin-Tits.\\
{\rm(i)}\qua $A_S^+$ est  simplifiable, c'est \`a dire que si $a,b,e_1,e_2$ sont dans $A_S^+$ et $ae_1b = ae_2b$ alors, $e_1 = e_2$.\\
{\rm(ii)}\qua Toute partie finie de $A_S^+$ poss\`ede un $\pgcd$ \`a gauche (et \`a droite).\\
{\rm(iii)}\qua Une partie finie de $A_S^+$ poss\`ede un $\ppcm$ \`a droite
 (resp.\ \`a gauche) si et seulement si elle poss\`ede un multiple commun \`a droite (resp.\ \`a gauche).\\
\end{Lemme}  
\subsection{Groupes d'Artin-Tits}
\begin{Lemme} [\cite{Cha1} th\'eor\`eme 2.6 et $\cite{Cha2}$ lemme 4.4]
{\relax}\ Soit $(A_S,S)$ un syst\`eme\break d'Artin-Tits de type sph\'erique. Si $g\in A_S$
 alors il existe $a,b$ dans $A_S^+$ uniques premiers entre eux \`a
 gauche ( not\'e $a\perp_\prec b$) tels que $g = a^{-1}b$. De plus si $c\in A_S^+$ tel que $cg\in A_S^+$ alors $c\succ a$. 
\end{Lemme} 
Nous appellerons  la d\'ecomposition $g = a^{-1}b$ l'\'ecriture
normale (\`a gauche) de $g$. On peut d\'efinir de la m\^eme fa\c con une
\'ecriture normale \`a droite.

\begin{Lemme}[\cite{VdL} th\'eor\`eme 4.13]\label{deftr}
Soit $(A_S,S)$ un syst\`eme d'Artin-Tits associ\'e \`a la matrice de
Coxeter $M = (m_{s,t})_{s,t\in S}$. Soit $T$ une partie de $S$ et $A_T$
le sous-groupe de $A_S$ engendr\'e par $T$. Alors $(A_T,T)$ est un
syst\`eme d'Artin-Tits associ\'e \`a la matrice $(m_{s,t})_{s,t\in T}$. De plus, si $T'$ est une autre partie de $S$, on a  $A_T\cap A_{T'} = A_{T\cap T'}$.\end{Lemme}

\subsection{lcm-homomorphismes}
\begin{Def}\label{deflcmcri}(\cite{Cri} d\'efinition 1.1)\qua Soit $A_S$ et $A_{S'}$ deux
groupes d'Artin-Tits. Si $\varphi: A_S\to A_{S'}$ est un homomorphisme tel que
$\varphi (A_S^+) \subset A^+_{S'}$, on note $\varphi^+$ la restriction de
$\varphi$ \`a $A^+_S$ et $A^+_{S'}$.  On dit que $\varphi^+$ (ou $\varphi$)
respecte les $\ppcm$ si\\(1)\qua $\forall s\in S$, $\varphi^+(s)\not= 1 $
et,\\(2)\qua $\forall s,t\in S$ on a $\varphi^+(s)\lor_\prec \varphi^+(t) = \varphi^+( s\lor_\prec t)$; \\avec la convention $\varphi(\infty) = \infty$.\end{Def}
Sous ces hypoth\`eses, on dira que $\varphi$ est un ``$\lcm$-homomorphisme''. 
\begin{Prop}[\cite{Cri} lemme 1.2 et Th\'eor\`eme 1.3]\label{Thcrisp}
Soit $\varphi: A_S\to A_{S'}$ un $\lcm$-homomorphisme et $U,V$ dans $A^+_S$;
alors $\varphi^+(U)\prec \varphi^+(V)$ si et seulement si $U\prec V$. En particulier $\varphi^+$ est injective. \end{Prop}
L'injectivit\'e peut aussi \^etre vu comme un cas particulier de la
proposition 5.4 de \cite{Deh1} sur les morphismes entre groupes munis de
pr\'esentations compl\'ement\'ees noeth\'eriennes et coh\'erentes. 

\begin{Que} Un $\lcm$-homomorphisme $\varphi:A_S\to A_{S'}$ est-il toujours
injectif?   \end{Que}
\section{LCM-homomorphismes}
Dans \cite{Cri} on d\'efinit une sous-famille de lcm-homomorphismes,
appel\'es $LCM$-homomorphismes; ceux-ci poss\`edent une r\'ealisation
g\'eom\'etrique naturelle. Cet\-te d\'efinition suppose que le graphe de $S$
soit sans liaisons infinies. Nous allons \'etendre
cette d\'efinition et montrer que le lemme clef 2.2 de \cite{Cri} est
encore vrai. La construction g\'eom\'etrique est identique \`a la construction de \cite{Cri}.\\ Commen\c cons par
rappeler la d\'efinition d'un LCM-homomorphisme.
\begin{Def}
Soit $(A_S,S)$ et $(A_{S'},S')$ deux syst\`emes d'Artin-Tits et soit $p$ une
 application de $S$ dans $ \mathcal{P}(S')-\{\emptyset\}$, les parties
 non vides de $S'$, telle
 que \\ {\bf(L0)}\qua si $s\not=t\in S$ alors $p(s)$ et $p(t)$ sont
 disjointes;\\ {\bf(L1)}\qua pour $s\in S$, $p(s)$ est de type sph\'erique;\\
{\bf (L2)}\qua  si $s\not=t\in S$ avec $m_{s,t}\not= \infty$, on a\\
 \begin{center} $[\Delta_{p(s)},\Delta_{p(t)}\rangle^{m_{s,t}} = [\Delta_{p(t)},\Delta_{p(s)}\rangle^{m_{s,t}} = \Delta_{p(s)}\lor_\prec \Delta_{p(t)}$ dans $A_{S'}^+$\end{center} 
{\bf(L3)}\qua si $s\not=t\in S$ avec $m_{s,t}= \infty$, alors
 \begin{center}$\left\{\begin {array}{l} \forall u\in p(s),\ \{u\}\cup p(t)\textrm{ n'est pas de type sph\'erique},\\ \forall u\in p(t),\ \{u\}\cup p(s)\textrm{ n'est pas de type sph\'erique.}\end{array}\right.$ \end{center}
 On peut alors d\'efinir un lcm-homomorphisme $\varphi_p$ de $A_S$ dans
 $A_{S'}$ par $\varphi_p(s)= \Delta_{p(s)}$ pour $s\in S$. Un morphisme provenant d'une telle construction sera appel\'e un LCM-homomorphisme.      
\end{Def} \`A la place de l'axiome {\bf(L3)}, on
 aurait pu se contenter, pour obtenir un lcm-homomorphisme, de l'axiome
 plus faible suivant: 

{\bf(L3$'$)}\qua si $s,t\in S$ et
 $m_{s,t}= \infty$, alors $p(s)\cup p(t)$ n'est pas de type
 sph\'erique

 mais l'axiome {\bf(L3)} a pour objectif de prouver la proposition
 \ref{lemQF}. \\Puisqu'un LCM-homomorphisme est un lcm-homomorphisme, on
 a le r\'esultat suivant:
\begin{Lemme} Soit $(A_S,S)$ et $(A_{S'},S')$ deux syst\`emes d'Artin-Tits et $\varphi_p$
un LCM-homomorphisme; alors $\varphi^+_p$ est injective de $A_S^+$ dans $A_{S'}^+$.\end{Lemme}  
Nous suivons dans la suite de cette partie le plan de \cite{Cri}. Nous
commen\c cons par rappeler deux lemmes techniques et par en d\'emontrer
un troisi\`eme; ceux-ci vont nous servir \'etablir la
proposition \ref{lemQF}. 
\begin{Lemme}\label{chaine}{\rm \cite{BrS,Cri2}}\qua
Soit $(A_S,S)$ un syst\`eme d'Artin-Tits et $T\subset S$. Soit $t\in T$, $w\in A_T^+$ et $x\in A_S^+$. Si $t\not\prec w$ mais $t\prec wx$ alors il existe $s\in T$ tel que $s\prec x$. 
\end{Lemme}
\begin{Lemme}[\cite{BrS} lemme 3.4 \label{lemred2}]Soit $(A_S,S)$ un
 syst\`eme d'Artin-Tits, soit
 $x\in A_{S,red}$ et $s\in S$; si $sx\not\in A_{S,red}$ alors $s\prec x$.  
\end{Lemme}
Ce lemme est une version du lemme d'\'echange.

\begin{Lemme}\label{lemred}Soit $(A_S,S)$ et $(A_{S'},S')$ deux syst\`emes
 d'Artin-Tits et $\varphi_p$ un LCM-homomorphisme. Soit $s,t\in S$ et $k\in \mathbb{N}$. Alors \begin{center} $[\varphi_p^+(s),\varphi_p^+(t)\rangle^k$ n'est pas r\'eduit $\Rightarrow m_{s,t}\not= \infty$ et $k > m_{s,t}$.  \end{center}
\end{Lemme} 

\begin{proof}[Preuve] Montrons  cette implication par contrapos\'ee.\\ Si $m_{s,t}\not=\infty$ et
$k\leq m_{s,t}$ alors $[\varphi_p^+(s),\varphi_p^+(t)\rangle^k$ divise
$[\varphi_p^+(s),\varphi_p^+(t)\rangle^{m_{s,t}}$; or $[\varphi_p^+(s),\varphi_p^+(t)\rangle^{m_{s,t}} = \Delta_{p(s)\cup p(t)}$ par l'axiome
\textbf{(L2)} et est ainsi r\'eduit ({cf.} la remarque qui suit le lemme \ref{lemsecpi}), donc $[\varphi_p^+(s),\varphi_p^+(t)\rangle^k$ est r\'eduit. Supposons maintenant $m_{s,t} = \infty$ et
montrons, par r\'ecurrence sur $k$, que
$[\varphi_p^+(s),\varphi_p^+(t)\rangle^k$  et $[\varphi_p^+(t),\varphi_p^+(s)\rangle^k$ sont r\'eduit
 pour tout $k$. Si $k= 0$ ou $k = 1$ c'est
clair; supposons donc $k\geq 2$ et que $[\varphi_p^+(s),\varphi_p^+(t)\rangle^k$
n'est pas r\'eduit. Par hypoth\`ese de r\'ecurrence
$[\varphi_p^+(t),\varphi_p^+(s)\rangle^{k-1}$ est  r\'eduit; donc il existe $C,D\in A^+_{p(s)}$
et $u\in p(s)$ tels que l'on a  $CuD = \varphi_p^+(s)$ avec
 $D[\varphi_p^+(t),\varphi_p^+(s)\rangle^{k-1}$  r\'eduit et
$uD[\varphi_p^+(t),\varphi_p^+(s)\rangle^{k-1}$ qui ne l'est pas. D'apr\`es le lemme
\ref{lemred2}, on a alors $u\prec D[\varphi_p^+(t),\varphi_p^+(s)\rangle^{k-1}$; d'autre part, $\varphi_p^+(s)$ est
r\'eduit, donc $uD$ l'est aussi et $u\not\prec D$. Par le lemme \ref{chaine}
appliqu\'e avec $T=p(s)$, $w = D$ et $x =[\varphi_p^+(t),\varphi_p^+(s)\rangle^{k-1} $, il existe
$v\in p(s)$ qui divise $[\varphi_p^+(t),\varphi_p^+(s)\rangle^{k-1}$ pour $\prec$. Mais ceci implique que $v\lor_\prec
\Delta_{p(t)}$ existe et donc que $\{v\}\cup p(t)$ est de type
 sph\'erique; ce qui
contredit l'axiome \textbf{(L3)}. Donc $[\varphi_p^+(s),\varphi_p^+(t)\rangle^k$ est
r\'eduit. Par sym\'etrie, on a aussi que $[\varphi_p^+(t),\varphi_p^+(s)\rangle^k$ est
r\'eduit.  \end{proof}             
\begin{Prop}\label{lemQF} Soit $A_S$ et $A_{S'}$ deux groupes
 d'Artin-Tits et $\varphi_p$ un LCM-homomorphisme. Alors la restriction de
 $\varphi^+_p$ \`a $A_{S,red}$ est un morphisme injectif
 dont l'image est incluse dans $A_{S',red}$.
\end{Prop}
Pour parler de cette propri\'et\'e de $\varphi_p$, on dira que $\varphi_p$ est
QF-injective.

\begin{proof}[Preuve] Puisque $\varphi_p$ est un LCM-homomorphisme, sa restriction $\varphi_p^+$ est
injective et la restriction \`a $A_{S,red}$ aussi. Il suffit donc de montrer que l'image
par $\varphi^+_p$ d'un \'el\'ement r\'eduit est r\'eduit. Soit $U\in A_{S,red}$
; on montre par r\'ecurrence sur la longueur de $U$ que $\varphi^+_p(U)$ est
r\'eduit. Si $ \ell(U) = 0$ ou $\ell(U) = 1$, alors le r\'esultat est
vrai. Supposons donc $\ell(U)\geq 2$. Dans ce cas, puisque $U$ est
r\'eduit, on peut \'ecrire $U = [s,t\rangle^mV$ avec $s,t\in S$ distincts et
$V\in A^+_S$ divisible pour $\prec$ ni par $s$ ni par $t$; de plus on a alors
 $2\leq m \leq
m_{s,t}$ et $V$ r\'eduit. On a alors que 
$[\varphi^+_p(s),\varphi^+_p(t)\rangle^m$  est r\'eduit d'apr\`es le lemme \ref{lemred}. Comme $V\in A_{S,red}$ et $\ell(V) < \ell(U)$, on a par
hypoth\`ese de r\'ecurrence $\varphi^+_p(V)$ r\'eduit. Supposons que
$\varphi^+_p(U)$ ne soit pas r\'eduit et proc\'edons comme dans la preuve
du lemme \ref{lemred}: il existe $C,D\in A^+_{p(s)\cup p(t)}$ et $u\in
p(s)\cup p(t)$ tel que $CuD = [\varphi^+_p(s),\varphi^+_p(t)\rangle^m$ et $D\varphi_p^+(V)$ est r\'eduit mais pas $uD\varphi_p^+(V)$; par
le lemme \ref{lemred2}, on a $u\prec D\varphi_p^+(V)$. D'autre part, comme
$[\varphi^+_p(s),\varphi^+_p(t)\rangle^m$ est r\'eduit, $uD$ l'est aussi et donc
$u\not\prec D$. Par
le lemme \ref{chaine}, il existe $v\in p(s)\cup p(t)$ qui divise
$\varphi_p^+(V)$ pour $\prec$. Supposons $v\in p(s)$; puisque
 $\Delta_{p(s)}\succ v$, $\varphi_p^+(sV)$ n'est pas r\'eduit. D'autre part,
$\ell(sV)<l(U)$, donc par hypoth\`ese de r\'ecurrence on a que $sV$ n'est pas
r\'eduit, ce qui implique par le lemme \ref{lemred2} que $s\prec V$; ceci est impossible par construction de $V$. Si $v\in
p(t)$ on proc\`ede de la m\^eme fa\c con pour obtenir une nouvelle
contradiction. Donc $U$ est r\'eduit.\end{proof}
\begin{Cor}\label{cor1}
Soit $\varphi_p: A_S\to A_{S'}$ un LCM-homomorphisme. Alors $\varphi_p$ induit un
 homomorphisme injectif $\varphi_{p,W}: W_S\to W_{S'}$.
\end{Cor}

\begin{proof}[Preuve] $\varphi_p(s^2) = \Delta_{p(s)}^2$ a pour image $1$ dans
$W_{S'}$; donc $\varphi_{p,W}$ est bien d\'efinie. D'autre part gr\^ace \`a la
proposition \ref{lemQF} et la section $Sec$, il est clair que $\varphi_{p,W}$
est injective. \end{proof}
\begin{Def} Soit  $\varphi:A_S\to A_{S'}$ un $\lcm$-homomorphisme. Pour
 $s\in S$ on pose $p_{\prec}(s) = \{t\in S' | t\prec \varphi(s)\}$ et
 $p_{\succ}(s) = \{t\in S' | \varphi(s)\succ t\}$. On dira que $\varphi$ est
 un $lcm$-homomorphisme sym\'etrique si $\forall s\in S$, on a
 $p_\prec(s) = p_\succ(s)$. Dans ce cas on notera simplement cet
 ensemble $p(s)$.
\end{Def}
Un LCM-homomorphisme est en particulier un lcm-homomorphisme sym\'etrique.
\begin{Lemme} \label{lemQF2}
Soit $\varphi: A_S\to A_{S'}$ un lcm-homomorphisme QF-injectif.
 Soit $U\in A_S^+$ r\'eduit et $s\in S$; s'il existe $t\in p_\succ(s)$ tel que
 $t\prec \varphi^+(U)$ alors $s\prec U$. Suposons de plus que $\varphi$ est sym\'etrique. Si $U$ et $V$ sont dans $A^+_S$ et r\'eduit, alors $\varphi^+(U)\land_\prec\varphi^+(V) = 1 \iff U\land_\prec V = 1$.   \end{Lemme}
 
\begin{proof}[Preuve]  Si $U$ est r\'eduit et
 $s\in S$ avec $s\not\prec U$, alors $sU$ est aussi r\'eduit. Donc
$\varphi^+(s)\varphi^+(U)$ est r\'eduit et aucun \'el\'ement de $p_\succ(s)$ ne
peut diviser $\varphi^+(U)$. Supposons maintenant que $\varphi$ est sym\'etrique et que $U$ et $V$ sont dans $A^+_S$ r\'eduit et diff\'erents de $1$. Par contreappos\'ee, il est clair que  $\varphi^+(U)\land_\prec\varphi^+(V) = 1 \Rightarrow U\land_\prec V = 1$. Montrons l'autre implication par l'absurde. Supposons que  $U\land_\prec V = 1$ mais $\varphi^+(U)\land_\prec\varphi^+(V) \neq 1$. Soit $t\in S'$ tel que $t\prec\varphi^+(U)\land_\prec\varphi^+(V)$  et notons $s\in S$  tel que $t\in p(s)$. En utilisant la premi\`ere partie du lemme, on a $t\prec \varphi^+(U)\Rightarrow s\prec U$ et $t\prec \varphi^+(V) \Rightarrow s\prec V$. Ce qui donne  finalement $s\prec U\land_\prec V$ et la contradiction voulue.  
\end{proof}
\begin{The} \label{cor2} Soit $\varphi: A_S\to A_{S'}$ un
 lcm-homomorphisme sym\'etrique QF-injectif. Alors
 $\varphi^+$ est compatible avec la forme normale: si $(g_1,\cdots,g_n)$ est la forme
 normale de $g\in A_S^+$ alors $(\varphi^+(g_1),\cdots,\varphi^+(g_n))$ est la
 forme normale de $\varphi^+(g)$. De plus, si $g_1,g_2\in A_S^+$ alors
 $$\left.\begin{array}{c l} (a) & \varphi^+(g_1)\lor_\prec \varphi^+(g_2) = \varphi^+(g_1\lor_\prec g_2)\\ (b) & \varphi^+(g_1)\land_\prec \varphi^+(g_2) = \varphi^+(g_1\land_\prec g_2).\end{array}\right.$$ 
\end{The}
Ce th\'eor\`eme est en particulier valable pour les LCM-homomorphismes
par la proposition \ref{lemQF}
et dans le cas des lcm-homomorphismes provenant des automorphismes de graphes
(\textit{cf.} th\'eor\`eme 14 de \cite{Cri2}). Le (a) est en fait connu
puisque Crisp a montr\'e dans le th\'eor\`eme 8 de \cite{Cri2}
que c'est d\'ej\`a le cas pour un lcm-homomorphisme.

\begin{proof}[Preuve] Puisque $\varphi$ est QF-injective, la suite
$(\varphi^+(g_1),\cdots,\varphi^+(g_n))$ a bien ses termes dans $A_{S',red}$. On montre le
r\'esultat par r\'ecurrence sur $n$; rappelons qu'au lemme
 \ref{lemalpha}, on a d\'efini une fonction $\alpha$ qui, \`a un \'el\'ement
 d'un mono\"\i de d'Artin-Tits associe son plus grand diviseur
 r\'eduit et qui v\'erifie que $\alpha(gh) = \alpha(g\alpha(h))$. Puisque $\alpha(\varphi^+(g_1)\cdots\varphi^+(g_n)) =
\alpha(\varphi^+(g_1)\alpha(\varphi^+(g_2)\cdots\varphi^+(g_n)))$, il suffit de
montrer le r\'esultat pour $n = 2$. Supposons donc $n=2$. Puisque $\varphi^+(g_1)$ est r\'eduit, il divise $\alpha(\varphi^+(g))$. Si
$\varphi^+(g_1)\not= \alpha(\varphi^+(g))$ alors il existe $t\in p(s)\subset S'$ avec
 $s\in S$ tel que $t\prec
\varphi^+(g_2)$ et $\varphi^+(g_1)\not\succ t$ mais dans ce cas $g_1\not\succ
s$ et par le lemme \ref{lemQF2}, on a $s\prec g_2$; ceci implique
que $g_1s$ est r\'eduit et divise $g$; ce qui contredit le fait que $g_1 =
\alpha(g)$. Donc $\varphi^+(g_1) = \alpha(\varphi^+(g))$.\\ Soit maintenant $g_1,g_2\in A^+_S$. Il est clair que
 $\varphi^+(g_1)\lor_\prec \varphi^+(g_2)\prec \varphi^+(g_1\lor_\prec g_2)$ et que
 $\varphi^+(g_1\land_\prec g_2)\prec \varphi^+(g_1)\land_\prec \varphi^+(g_2)$. D'autre part, si on note $g_1 = (g_1\land_\prec g_2)h_1$ et $g_2 = (g_1\land_\prec g_2)h_2$ alors $h_1\land_\prec h_2 = 1$. Il est facile de voir que $h_1\land_\prec h_2 = 1 \iff \alpha(h_1) \land_\prec \alpha(h_2) = 1$ et que  $\varphi^+(h_1)\land_\prec \varphi^+(h_2) = 1 \iff \alpha(\varphi^+(h_1)) \land_\prec \alpha(\varphi^+(h_2)) = 1$ (voir le lemme \ref{lemalpha}). On en d\'eduit alors par le lemme \ref{lemQF2} et la premi\`ere partie du th\'eor\`eme  que $\varphi^+(h_1)\land_\prec \varphi^+(h_2) = 1$, ce qui montre le (b).\\ Le (a) se montre de la m\^eme fa\c con : si on \'ecrit $g_1\lor_\prec g_2 = g_1h_1 = g_2h_2$, on a $h_1\land_\succ h_2 = 1$ et $\varphi^+(h_1)\land_\succ\varphi^+(h_2) = 1.$ \end{proof}  

\noindent Remarquons que l'on a vraiment eu besoin du fait que le lcm-homomorphisme est
sym\'etrique et que l'on ne
peut esp\'erer \'etendre ce r\'esultat \`a tous les lcm-homomorphismes QF-injectifs
comme le montre l'exemple suivant: soit $\varphi$ le lcm-homomorphisme du
groupe libre \`a un g\'en\'erateur $\langle t \rangle$ dans le groupe libre \`a deux
g\'en\'erateurs $\langle x,y \rangle$ qui envoie $t$ sur $xy$. Alors $\varphi(t^2)= xyxy$
est r\'eduit  donc $\alpha(xyxy) = xyxy$ mais $\varphi(\alpha(t^2)) =
\varphi(t) = xy$. 
\begin{Cor} \label{cor3} Soit $\varphi_p: A_S\to A_{S'}$ un LCM-homomorphisme entre
 groupes d'Artin-Tits de type sph\'erique. Alors
 $\varphi_p^+$ conserve l'\'ecriture normale: si $g_1^{-1}g_2$ est l'\'ecriture
 normale de $g\in A_S$ alors $\varphi^+_p(g_1)^{-1}\varphi^+_p(g_2)$ est l'\'ecriture normale de $\varphi_p(g)$. 
\end{Cor}

\begin{proof}[Preuve] C'est clair par la seconde partie du th\'eor\`eme \ref{cor2}. \end{proof}
\subsection{Une r\'ealisation g\'eom\'etrique}
Le r\'ealisation g\'eom\'etrique se construit maintenant exactement comme dans
\cite{Cri}; on se contente donc d'introduire les notations utiles
et  d'\'enoncer le r\'esultat.\\       
Soit $(A_S,S)$ un syst\`eme d'Artin-Tits, $(W_S,S)$ son syst\`eme de Coxeter
associ\'e et posons $\mathcal{S}_{f,S} =
\{T\subset S | W_T$ est fini $\}$. On ordonne l'ensemble $W_S\times
\mathcal{S}_{f,S}$ par la relation: $$(w_1,T_1)\leq (w_2,T_2)\textrm{ si
}\left\{\begin{array}{l} w_1W_{T_1}\subset
w_2W_{T_2}\\ \ell(\Delta_{T_1})+\ell(w_1^{-1}w_2) =
\ell(\Delta_{T_1}w_1^{-1}w_2) \end{array}\right.$$ et on note $\widetilde{\Sigma}_S$ la
r\'ealisation g\'eom\'etrique du complexe d\'eriv\'e de cet
ensemble ordonn\'e; c'est un complexe simplicial.  Le complexe
$\widetilde{\Sigma}_S$ s'appelle le complexe de Salvetti et a \'et\'e introduit dans \cite{Sal}. Le groupe $W_S$ agit
simplicialement et librement sur $\widetilde{\Sigma}_S$ (par multiplication
\`a gauche sur le premier facteur); on note $Z_S = \widetilde{\Sigma}_S/W_S$
l'espace des orbites pour cette action. Il est connu que $\pi_1(Z_S)\cong
A_S$.
\begin{Lemme} [Proposition 3.2 de \cite{Cri}]\label{lemPhi} Soit $A_S$ et $A_{S'}$
 deux groupes d'Artin-Tits. Soit $\varphi_p$ un LCM-homomorphisme. Alors
 l'application $\widetilde{p}: W_S \times \mathcal{S}_{f,S}\to W_{S'}\times \mathcal{S}_{f,S'}$ d\'efinie par $\widetilde{p}(w,T) = (\varphi_{p,W}(w),p(T))$ pour $(w,T)\in W_S\times \mathcal{S}_{f,S}$ pr\'eserve strictement l'ordre et induit une application $\widetilde{\Phi}$ qui est simpliciale injective et $(W_S,W_{S'})$-\'equivariante de $\widetilde{\Sigma}_S$ dans $\widetilde{\Sigma}_{S'}$.
\end{Lemme}
\begin{Prop}[Th\'eor\`eme 3.4 de \cite{Cri}]\label{prorep}
Soit $\varphi_p:A_S\to A_{S'}$ un LCM-hom\-o\-morphisme; identifions $A_S,A_{S'}$
 respectivement \`a $\pi_1(Z_S)$ et $\pi_1(Z_{S'})$. Alors l'application
 $\widetilde\Phi $ du lemme \ref{lemPhi} induit une
 application simpliciale injective $\Phi:Z_S\to Z_{S'}$ telle que
$$\varphi_p = \Phi_*,$$ o\`u $\Phi_*$ d\'esigne l'application induite par $\Phi$ entre  les groupes fondamentaux.\end{Prop}
 \section{Injectivit\'e des LCM-homomorphismes}
Nous commen\c cons par introduire le complexe de Deligne puis la notion
de complexe de cubes. C'est gr\^ace \`a ces objets que nous allons
\'etablir l'injectivit\'e des LCM-homomorphismes pour les
groupes d'Artin-Tits de type FC.   
\subsection{Le complexe de Deligne et les espaces de cubes}
Le complexe de Deligne a \'et\'e d\'efini pour la premi\`ere fois par
 Deligne dans \cite{Del} pour le cas des groupes d'Artin-Tits de type sph\'erique. Cette construction a
 \'et\'e ensuite g\'en\'eralis\'ee par Charney et Davis dans \cite{ChD1} et a
 permis de montrer plusieurs r\'esultats sur les groupes d'Artin-Tits (par
 exemple dans \cite{Cha1,Cha2,ChD2,ChD1,Cri2}).\\
 Soit $(A_S,S)$ un syst\`eme
 d'Artin-Tits. Rappelons  que $$\mathcal{S}_{f,S} = \{T\subset S;\ A_T\textrm{ est
 de type sph\'erique}\}$$ et posons $$A_S\mathcal{S}_{f,S} = \{xA_T;\ x\in
 A_S\textrm{ et }T\in \mathcal{S}_{f,S}\}.$$

Rappelons que si $(P,<)$ est un ensemble partiellement ordonn\'e, son complexe de drapeaux est le complexe simplicial abstrait qui a pour sommets les \'el\'ements de $P$ et pour simplexes les suites finies croissantes d'\'el\'ements de $P$. Un sous-simplexe d'un simplexe, donc d'une suite, est alors une sous-suite de celle-ci.\\
 On appelle ``complexe de Deligne'' le
 complexe de drapeaux $D_S$ obtenu \`a partir de l'ensemble $A_S\mathcal{S}_{f,S}$ munie de
 l'inclusion comme ordre partiel (notons que $xA_X\subset yA_Y \iff X\subset Y$ et
 $y^{-1}x\in A_Y$). C'est un complexe simplicial abstrait. On peut
 associer \`a ce complexe une r\'ealisation g\'eom\'etrique sous forme
 de complexe simplicial euclidien par morceaux. On identifie souvent le
 complexe et sa r\'ealisation, bien que celle-ci ne soit pas unique. Le groupe $A_S$ agit sur son complexe de Deligne
 par multiplication \`a gauche. Si l'on note $K_S$, le sous-complexe de
 drapeaux associ\'e \`a $\mathcal{S}_{f,S}$, alors $K_S$ est un domaine
 fondamental pour l'action de $A_S$ sur $D_S$. Si l'on munit $K_S$ d'une r\'ealisation g\'eom\'etrique, on peut \'etendre celle-ci \`a $D_S$ tout entier via l'action de $A_S$. Dans ce cas, $A_S$ agit naturellement par isom\'etries sur la r\'ealisation g\'eom\'etrique de $D_S$. \\

\noindent  Dans \cite{ChD1}, on associe au complexe de Deligne une r\'ealisation g\'eom\'etrique
 li\'ee \`a la m\'etrique dite ``de
 Moussong''; nous ne d\'etaillons pas ici cette construction. Disons
 simplement que cette r\'ealisation est conjecturalement la plus
 adapt\'ee pour tous les groupes d'Artin-Tits. Lorsque $A_S$ est de type FC (\textit{cf.} d\'efinition \ref{deffc}), on associe \`a $D_S$ une r\'ealisation
 g\'eom\'etrique cubique gr\^ace \`a $K_S$, comme vu ci-dessus. La structure cubique de $K_S$ est donn\'ee par les sous-complexes associ\'es aux
sous-groupes paraboliques standards $A_X$ de type sph\'erique dont les r\'ealisations g\'eom\'etriques sont alors des cubes de $\mathbb{R}^n$ pour $n = |X|$.

\noindent Un complexe de cubes est un complexe poly\'edrique o\`u les faces
ferm\'ees sont
des cubes d'un espace euclidien; celles-ci sont appel\'ees les cubes
du complexe. Si $(A_S,S)$ est un syst\`eme d'Artin-Tits
de type $FC$ alors la r\'ealisation cubique de son complexe de Deligne est naturellement munie d'une
structure de complexe de cubes sous-jacente qui consiste \`a ne garder
que les ar\^etes $\{aA_X,aA_Y\}$ telles que $aA_X\subset aA_Y$ avec $a\in A_S$ et
$Y = X\cup\{s\}$ dans $\mathcal{S}_{f,S}$ avec $s\in S-X$. Si $K$ est un
cube de dimension $n$, alors il poss\`ede pour l'inclusion un sommet minimal $aA_R$ et un sommet maximal
$aA_T$, avec $a\in A_S$, et avec $R\subset T$ dans
$\mathcal{S}_{f,S}$ tels que $\#(T-R) = n$. L'ensemble des sommets de $K$ est alors $\{aA_X|
R\subset X\subset T\}$ et il forme un treillis pour l'inclusion. On notera $K = K(aA_R,aA_T)$.    
\begin{Def} \label{defchepol} Soit $D$ un complexe de cubes.\\ {\rm(i)}\qua Si
 $C_1,C_2,\cdots, C_n$
sont des cubes de $D$, on note $span(C_1,C_2,\cdots,C_n)$ le plus petit cube,
s'il existe, qui contient les $C_i$. On dit que c'est le cube
tendu par les $C_i$.\\ {\rm(ii)}\qua Soit $C$ un cube de $D$. On
 appelle \'etoile de $C$, et on note $Et(C)$, le sous-complexe $Et(C)= \bigcup_{C\subset C'} C'$.\\  
{\rm(iii)}\qua Soit $x,y$ deux sommets de $D$. On dit que la suite $x_1,\cdots, x_n$  de sommets de $D$ est un chemin cubique normal de $x$ \`a $y$ si $$\left\{\begin{array}{c l} (a) & x =x_0; y= x_n ;\\ (b) & \forall i\in\{1,\cdots, n\}, C_i = span(x_{i-1},x_i) \textrm{ existe};\\ (c) & \forall i\in \{0,\cdots,n-1\}, Et(C_i)\cap C_{i+1} = \{x_i\}.\end{array}\right.$$\end{Def}
\begin{Prop}[\cite{NiR}, \cite{ChA} th\'eor\`eme 4.4] Soit $(A_S,S)$ un syst\`eme d'Artin-Tits de type FC; on munit $D_S$ de sa structure de complexe de cubes. Soit $xA_X,yA_Y$ deux
 sommets de $D_S$; alors il existe un unique chemin cubique normal de $xA_X$
 \`a $yA_Y$.\end{Prop}
$$\epsfysize 4.5cm\epsfbox{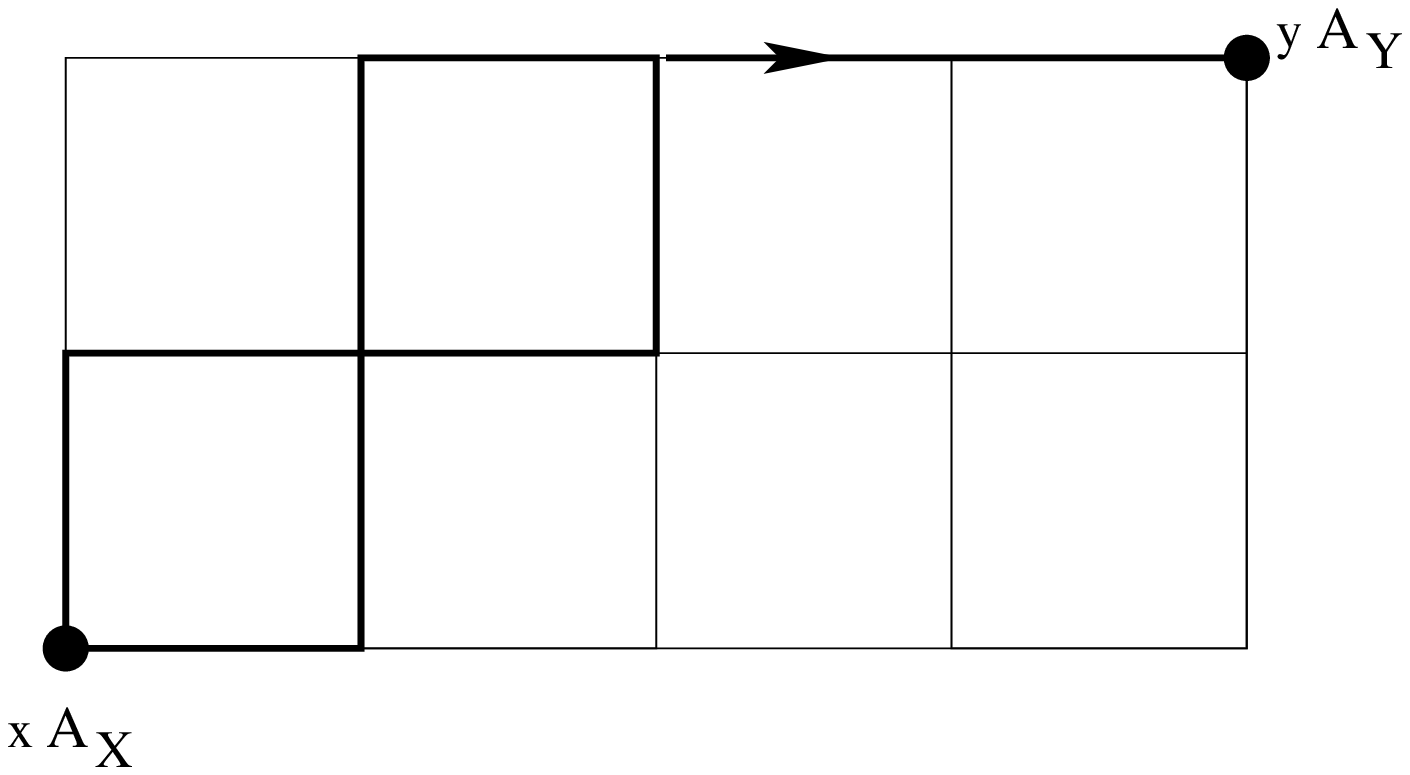}$$
\medskip
\centerline{\small Chemin cubique normal de $xA_X$ \`a $yA_Y$}

\noindent Cet \'enonc\'e est li\'e au fait que la r\'ealisation g\'eom\'etrique
 cubique d'un groupe d'Artin-Tits de type FC est  un espace
 CAT(0) (\textit{cf.} \cite{ChD1}). Nous ne d\'etaillons pas ici cette
 notion qui ne nous servira pas directement. On pourra se r\'ef\'erer \`a
 \cite{BrH} pour plus d'explications sur les espaces m\'etriques \`a courbure n\'egative.\\

\noindent Le lemme suivant permet de mieux comprendre la structure cubique du complexe de Deligne.   
\begin{Lemme}[\cite{ChA} lemme 4.1]\label{lemspa} Soit $(A_S,S)$ un
 syst\`eme d'Artin-Tits  de type FC et $D_S$ son complexe de Deligne munit de sa structure de
 complexe de cubes. Soit
 $K_1 = K(aA_{R_1},aA_{T_1})$ et $K_2 = K(bA_{R_2},bA_{T_2})$ deux cubes de
 $D_S$. Alors, $$span(K_1,K_2)\textrm{ existe }\iff T_1\cup T_2 \in \mathcal{S}_{f,S} \textrm{ et }aA_{R_1}\cap bA_{R_2}\neq \emptyset.$$     
De plus, dans ce cas, $$span(K_1,K_2) = K(c A_{R_1\cap R_2},c A_{T_1\cup T_2})$$ avec $c\in aA_{R_1}\cap bA_{R_2}$.\end{Lemme}
\subsection{Injectivit\'e des  LCM-homomorphismes}
Dans cette partie, on se donne $(A_S,S)$ et $(A_{S'},S')$
 deux syst\`emes d'Artin-Tits de type FC et $\varphi_p : A_S\to A_{S'}$ un LCM-homomorphisme. On identifie les complexes de Deligne $D_S$ et $D_{S'}$ avec
 leurs r\'ealisations g\'eom\'etriques cubiques. On d\'efinit
 une application simpliciale continue $\Phi_p:D_S\to D_{S'}$ en
 posant  $\Phi_p(xA_X) = \varphi_p(x)A_{p(X)}$. Cette application est bien
 d\'efinie car si $X\subset Y\subset S$ alors $p(X)\subset p(Y)\subset S'$.
L'\'enonc\'e du lemme suivant est largement inspir\'e du lemme 18 de \cite{Cri2}.   
\begin{Lemme}\label{lemccri} On a: \\{\rm(i)}\qua la restriction
 $\Phi_p|_{K_S}$ de $\Phi_p$ \`a $K_S$ est injective et son image
 est dans $K_{S'}$.\\{\rm(ii)}\qua Si $\Phi_p$
 est injective alors $\varphi_p$ est injective.\\{\rm(iii)}\qua $Im \varphi_p = stab(Im
 \Phi_p)$ et $Im \Phi_p = Im \varphi_p \cdot \Phi_p(K_S)$.
\end{Lemme}
\begin{Que} Il n'est pas tr\`es compliqu\'e de voir que la d\'efinition
 de $\Phi_p$ et le lemme \ref{lemccri} sont ind\'ependants du fait
 que les groupes d'Artin-Tits sont de type FC (c'est ind\'ependant de la
 r\'ealisation g\'eom\'etrique). Une question naturelle est donc de
 savoir si la r\'eciproque du (ii) est vraie en g\'en\'eral. 
\end{Que}

\begin{proof}[Preuve] {\rm(i)}\qua La restriction de $\Phi_p$ \`a $K_S$ est injective sur les sommets puisque $p$
envoie des g\'en\'erateurs distincts sur des parties disjointes. Par
simplicialit\'e, elle est donc injective sur $K_S$ et par d\'efinition
$\Phi_p(K_S)\subset K_{S'}$. Pour montrer le (ii), il suffit de
 consid\'erer l'orbite de $A_\emptyset$. \\ 
{\rm(iii)}\qua L'inclusion $Im \varphi_p \subset stab(Im
 \Phi_p)$ et l'\'egalit\'e  $Im \Phi_p = Im \phi_p \cdot \Phi_p(K_S)$ sont claires. D'autre part, on a $1A_\emptyset \in Im(\Phi_p)$, donc si $w \in stab(Im(\Phi_p))$
alors $w(1A_\emptyset) = wA_{\emptyset}\in Im(\Phi_p)$ et $w\in Im \varphi_p.$    \end{proof}
Avant d'\'enoncer la proposition importante \ref{Proccn} nous commen\c
cons par un lemme utile pour la preuve de celle-ci.
\begin{Lemme}\label{lemcons} Soit $A_S$ un groupe d'Artin-Tits.\\
{\rm(i)}\qua  Soit $w\in A_S^+$,  $X\subset S$, et $s\in S$. Si $w\in A_X^+$ et $s$ appara\^\i t dans
 une \'ecriture de $w$ alors $s\in X$.\\{\rm(ii)}\qua
 Soit $\varphi_p:A_S\to A_{S'}$ un LCM-homomorphisme. Soit
$R\subset S$ et $Y\subset S'$. On pose $Z = \{s\in S| p(s)\subset Y\}$. Soit $w\in A^+_S$ alors $$\exists \alpha\in
 A^+_{p(R)}\ ,\exists \beta\in A^+_{Y}\ ,\  \varphi_p(w) = \alpha\beta \Rightarrow \exists u\in A^+_{R}\ ,\exists\  v\in A^+_{Z}\ ,\  w = uv.$$ En particulier  si $R = \emptyset$ ( donc $\alpha = 1$ et $\varphi_p(w)\in A^+_{Y}$) alors $w\in A^+_{Z}$.   \end{Lemme}
Au cours de la preuve de ce lemme et de la proposition suivante, nous
 aurons besoin des notions d'\'el\'ement $X$-r\'eduit et r\'eduit-$X$:
\begin{Def} Soit $A_S$ un groupe d'Artin-Tits; soit $X$ est une partie
de $S$ et $w\in A_S^+$. On dit que $w$ est $X$-r\'eduit (resp.\ r\'eduit-$X$) s'il n'est divisible pour $\prec$ (resp.\ $\succ$) par aucun \'el\'ement de $X$. \end{Def} 

\begin{proof}[Preuve du lemme \ref{lemcons}]{\rm(i)}\qua Supposons que $w\in A_X^+$. Cela signifie qu'il existe un
 repr\'esentant de $w$ dans $A_S^+$ o\`u n'apparaissent que des \'el\'ements de $X$; mais  si $s$ appara\^\i t dans une \'ecriture, alors il appara\^\i t
dans toutes car les relations de tresses ne font ni appara\^\i tre ni
dispara\^\i tre de g\'en\'erateurs. Donc $s\in X$.\\ {\rm(ii)}\qua Supposons
$l(\alpha) = 0$ pour commencer et
remarquons que si $X,X'\subset S$ alors $X\subset X'\iff p(X)\subset p(X')$. Soit $X$ minimal tel que $w\in A^+_{X}$: $w = s_1\cdots s_n$ avec $X = \{s_1,\cdots s_n\}$ (les $s_i$ ne sont pas forc\'ement distincts). Alors $\varphi^+_p(w) = \Delta_{p(s_1)}\cdots \Delta_{p(s_n)}$. D'o\`u $p(X) = \bigcup_{i=1}^{n}{p(s_i)}\subset Y$ par le (i) et donc $X\subset Z$ par d\'efinition de $Z$.\\ Revenons maintenant au cas g\'en\'eral : $\varphi_p(w) = \alpha\beta$.\\  Quitte \`a
 modifier $\alpha$ et $\beta$, on peut supposer que $\beta$ est
 $p(R)$-r\'eduit. On proc\`ede par r\'ecurrence sur $l(\alpha)$. Si
 $l(\alpha) = 0$, le r\'esultat est vrai par le d\'ebut de la
 preuve. Supposons donc  $l(\alpha)\geq 1$. Soit $t\in p(R)$
tel que $t\prec \alpha$ et $s\in R$ tel que $t\in p(s)$. Alors, $t\prec \varphi_p(w)$ et
comme $\varphi_p$ conserve la forme normale, $s\prec w$. Donc
$\Delta_{p(s)}\prec \alpha\beta$; mais comme $p(s)\subset p(R)$ et
$\beta$ est $p(R)$-r\'eduit, on a alors $\Delta_{p(s)}\prec \alpha$ car tout
\'el\'ement de $p(s)$ divise $\alpha$ par le lemme \ref{chaine}. On peut donc simplifier par $s$ dans $w$ et par $\Delta_{p(s)}$ dans $\alpha$ pour appliquer l'hypoth\`ese de r\'ecurrence.  \end{proof}    

\begin{Prop}\label{Proccn}
Soit $xA_X$ et $yA_Y$ deux sommets de $D_S$. Alors $\Phi_p$ envoie le chemin
 cubique normal de $xA_X$ \`a $yA_Y$ sur le chemin cubique normal de
 $\Phi_p(xA_X)$ \`a $\Phi_p(yA_Y)$.
\end{Prop}

\begin{proof}[Preuve] Soit $xA_X$ et $yA_Y$ des sommets de $D_S$. Notons
$x_0A_{X_0} = xA_X, x_1A_{X_1}$, $\cdots, x_nA_{X_n} = yA_Y$ l'unique chemin
cubique normal de $xA_X$ \`a $yA_Y$ dans $D_S$. Posons $R_i = X_{i-1}\cap X_i$
 et $T_i = X_{i-1}\cup X_i$. Soit $a_i\in x_{i-1}A_{X_{i-1}} \cap x_iA_{X_i}$; on a $$K(a_iA_{R_i},a_iA_{T_i}) = span (x_{i-1}A_{X_{i-1}},x_iA_{X_i}).$$ L'image par $\Phi_p$ du chemin cubique normal de $x$ \`a $y$ est la suite de sommets de $D_{S'}$: $\varphi_p(x_0)A_{p(X_0)} = \varphi_p(x)A_{p(X)},\ \varphi_p(x_1)A_{p(X_1)}$, $\cdots$, $\varphi_p(x_n)A_{p(X_n)} = \varphi_p(y)A_{p(Y)}$. Posons $C_i = span(\varphi_p(x_{i-1})A_{p(X_{i-1})},\varphi_p(x_i)A_{p(X_i)})$. Celle-ci existe et est en fait clairement \'egale \`a $C_i = K(\varphi_p(a_i)A_{p(R_i)},\varphi_p(a_i)A_{p(T_i)}).$\\  Nous devons prouver que cette suite de sommets v\'erifie les axiomes (a),(b), et (c) de la d\'efinition \ref{defchepol} (iii). Le (a) et le (b) sont triviaux; montrons donc le (c) par l'absurde. Soit $i$ tel que $Et(C_i)\cap C_{i+1}\neq \{\varphi_p(x_i)A_{p(X_i)}\}$ et soit $\varphi_p(a_{i+1})A_{Y}$ un sommet de $Et(C_i)\cap C_{i+1}$ distinct de ${\varphi_p(x_i)A_{p(X_i)}}$. Par le lemme \ref{lemspa}, on a $$span(\varphi_p(x_i)A_{p(X_i)},\varphi_p(a_{i+1})A_{Y}) = K(\varphi_p(a_{i+1})A_{p(X_i)\cap Y},\varphi_p(a_{i+1})A_{p(X_i)\cup Y}).$$ D'autre part, $p(X_i)\cup Y\neq p(X_i)$ ou $p(X_i)\cap Y\neq p(X_i)$; quitte \`a remplacer $Y$ par $p(X_i)\cap Y$ ou $p(X_i)\cup Y$, on peut donc se ramener soit au cas o\`u $\varphi_p(x_i)A_{p(X_i)}\varsubsetneq\ \varphi_p(a_{i+1})A_Y$ soit au cas o\`u $\varphi_p(a_{i+1})A_Y\varsubsetneq \varphi_p(x_i)A_{p(X_i)}$. \\

\noindent Premier cas: supposons que $\varphi_p(a_{i+1})A_Y\varsubsetneq \varphi_p(x_i)A_{p(X_i)}$ et posons $Z = \{s\in S; p(s)\subset  Y\}$. On a
 donc $Z\varsubsetneq X_i$.\\
Puisque $a_iA_{X_i} = a_{i+1}A_{X_i} = x_iA_{X_i}$, on a
 $a^{-1}_{i+1}a_{i}\in A_{X_i}$: $a^{-1}_{i+1}a_{i} =u_1^{-1}u^{-1}vv_1$
                                            avec $u,v,u_1,v_1\in A_{X_i}^+$ et $\left\{\begin{array}{l}uu_1\textrm{ et }vv_1\textrm{ premiers entre eux pour }\prec, \\v_1\in A^+_{R_i}\textrm{ et }v\textrm{ r\'eduit-}R_i,\\ u_1\in A^+_{Z}\textrm{ et }u\textrm{ r\'eduit-}Z.\end{array}\right.$\\ Montrons que $a_iA_{R_i}\cap a_{i+1}A_{Z}\neq\emptyset$; ceci est \'equivalent \`a montrer que $vA_{R_i}\cap uA_{Z}\neq\emptyset$. Or, par hypoth\`ese $span(\varphi_p(a_{i})A_{p(R_i)}, \varphi_p(a_{i+1})A_Y)$ existe et par le lemme \ref{lemspa} on en d\'eduit que $\varphi_p(a_{i})A_{p(R_i)} \cap \varphi_p(a_{i+1})A_Y \neq \emptyset$;
 c'est-\`a-dire qu'il existe $\alpha\in A_{p(R_i)}$ et $\beta\in A_Y$
 tel que $\varphi_p(v)\varphi_p(v_1)\alpha = \varphi_p(u)\varphi_p(u_1)\beta$. On peut \'ecrire
 $\varphi_p(v_1)\alpha = \alpha_1\alpha_2^{-1}$ avec
 $\alpha_1,\alpha_2\in A_{p(R_{i})}^+$  et $\alpha_1\land_\succ\alpha_2 = 1$; de m\^eme, on peut \'ecrire $\varphi_p(u_1)\beta = \beta_1\beta_2^{-1}$ avec $\beta_1,\beta_2\in A_{Y}^+$ et $\beta_1\land_\succ\beta_2 = 1$. Ceci nous donne $\varphi_p(v) \alpha_1\alpha_2^{-1} = \varphi_p(u)\beta_1\beta_2^{-1}$; mais l'\'ecriture $\varphi_p(v) \alpha_1\alpha_2^{-1}$  est normale: $\varphi_p(v) \alpha_1$ et $\alpha_2$ sont premiers entre eux \`a droite car $\alpha_1$ et $\alpha_2$ le sont, $\varphi_p(v)$ est r\'eduit-$p(R_i)$ car $v$ est r\'eduit-$R_i$ ({\it cf.} le lemme \ref{lemQF2}), et on utilise l'analogue \`a droite du lemme \ref{chaine}.
 Donc il existe $w$, a priori dans $A^+_{p(X_i)}$, tel que $\alpha_2w = \beta_2$ et $\varphi_p(v)\alpha_1w = \varphi_p(u)\beta_1$. L'\'egalit\'e $\alpha_2w = \beta_2$ impose $w\in A^+_Y$. Maintenant, en simplifiant $\alpha_1w$ et $\beta_1$ par leur pgcd \` a droite et en prenant pour $a$ le repr\'esantant r\'eduit-$Y$ de $\alpha_1A_Y^+$, on d\'eduit de l'\'egalit\'e  $\varphi_p(v)\alpha_1w = \varphi_p(u)\beta_1$ qu'il existe $a\in A^+_{p(R_i)}$ et $b,c\in A^+_{Y}$ tels que  que $a$ est r\'eduit-$Y$ et $ac\land_\succ b = 1$ avec $\varphi_p(v)ac = \varphi_p(u) b $. Cette \'egalit\'e implique que $b =\varphi_p(u)^{-1}(\varphi_p(u)\lor_\prec \varphi_p(v))z  = \varphi_p(u^{-1}(u\lor_\prec v))z $  et $ac = \varphi_p(v)^{-1}(\varphi_p (u)\lor_\prec\varphi_p (v))z = \varphi_p(v^{-1}(u\lor_\prec v))z$ avec $z\in A^+_Y$; mais puisque $ac\land_\succ b = 1$, on a $z = 1$. On peut maintenant appliquer le lemme \ref{lemcons}(ii) : $u^{-1}(u\lor_\prec v)\in A^+_{Z}$ et il existe $x\in A^+_{R_i},\ y\in A^+_Z$ tel que $v^{-1}(u\lor_\prec v) = xy$.  D'o\`u $v A_{R_i}\cap u A_{Z}\neq \emptyset$. De plus, $Z\cup T_{i} = T_{i}\in \mathcal{S}_{f,S}$, donc par la proposition \ref{lemspa}, on a que $span(K(a_iA_{R_i},a_iA_{T_i}),a_{i+1}A_{Z})$ existe; si l'on montre que  que l'on a aussi $a_{i+1}A_Z\in K_{i+1}$ on aura alors une contradiction avec le fait que la suite des $x_iA_{X_i}$ est un chemin
 cubique normal, puisque $Z$ et $X_i$ sont distincts. Montrons ce dernier point. Ceci revient \`a voir que $R_{i+1}\subset Z\subset T_{i+1}$. Mais d'une part on a $Z\subset X_i\subset T_{i+1}$, et d'autre part on a $p(R_{i+1})\subset Y$ puisque $\varphi_p(a_{i+1})A_Y$ est un sommet de $C_{i+1}$. Par definition de $Z$, cela implique $R_{i+1}\subset Z$.\\
Deuxi\`eme cas: supposons que
 $\varphi_p(x_i)A_{p(X_i)}\varsubsetneq \varphi_p(a_{i+1})A_Y$ (et
 distincts) et posons maintenant $Z = \{s\in S; p(s)\cap Y\neq \emptyset\}$. On
 a alors $R_i\subset X_i \varsubsetneq Z$ et $a_{i}A_{R_i} \cap a_{i+1}A_Z = a_{i}A_{R_i} \neq \emptyset$. De plus, puisque $span(\varphi_p(a_{i})A_{p(T_{i})},\varphi_p(a_{i+1})A_Y)$ existe, on a que $p(T_i)\cup Y$ est de type sph\'erique (en particulier sans liaison infinie);
 ceci implique par l'axiome \textbf{(L3)} que $T_i\cup Z$
 n'a pas non plus de liaison infinie. Il est donc de type sph\'erique
 puisque $S$ est de type $FC$. Comme dans le premier cas, on en d\'eduit
 par le lemme \ref{lemspa}, que $span(K(a_iA_{R_i},a_iA_{T_i}), a_{i+1}A_Z)$ existe.  Comme dans le premier cas, il nous reste \`a voir que $a_{i+1}A_Z\in K_{i+1}$, c'est \`a dire $R_{i+1}\subset Z\subset T_{i+1}$, pour obtenir une nouvelle contradiction puisque $Z$ et $X_i$ sont encore une fois distincts. Tout d'abord, on a $R_{i+1} \subset X_i\subset Z$. Ensuite puisque $\varphi_p(a_{i+1})A_Y$ est un sommet de $C_{i+1}$, on a $Y \subset p(T_{i+1})$; ce qui implique $p(Z)\subset p(T_{i+1})$ et finalement $Z\subset T_{i+1}$ par l'axiome {\bf(L0)} de la d\'efinition \ref{defLcm}.      
\end{proof}
$$\epsfysize 5.7cm\epsfbox{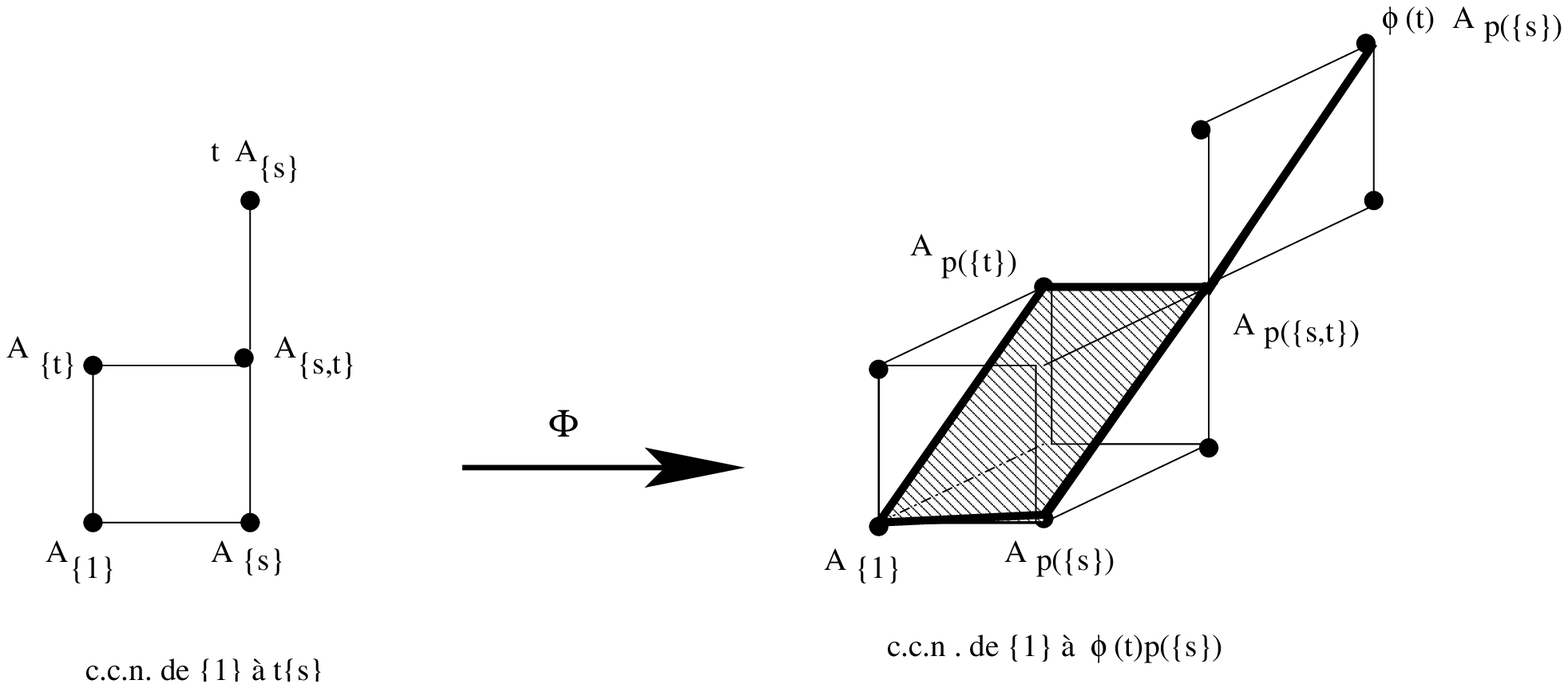}$$
\centerline{\small Un exemple}
\begin{Cor}\label{corinj} Soit $(A_S,S)$ et $(A_{S'},S')$ deux syst\`emes d'Artin-Tits de type FC
et $\varphi_p: A_S\to A_{S'}$ un LCM-homomorphisme. Alors, $\Phi_p$ est injective.
\end{Cor}

\begin{proof}[Preuve] Par la proposition pr\'ec\'edente, $\Phi_p$ est
         injective sur les sommets et est donc injective.\end{proof}
\begin{The} Soit $(A_S,S)$ et $(A_{S'},S')$ deux syst\`emes d'Artin-Tits de type FC
et $\varphi_p: A_S\to A_{S'}$ un LCM-homomorphisme. Alors $\varphi_p$ est injective.
\end{The}

\begin{proof}[Preuve] On applique le corollaire \ref{corinj} et le
         lemme \ref{lemccri}(ii).\end{proof}

\Addresses
\recd

\end{document}